\documentclass[11pt]{article}
\usepackage[margin=1in]{geometry}
\usepackage{amsmath,amssymb,amsthm,mathtools}
\usepackage{mathrsfs} 
\usepackage{esint} 
\usepackage[T1]{fontenc}
\usepackage{hyperref}
\newtheorem{theorem}{Theorem}[section]
\newtheorem{lemma}[theorem]{Lemma}
\newtheorem{proposition}[theorem]{Proposition}

\newtheorem{definition}{Definition}[section] 
\usepackage{enumitem}
\theoremstyle{remark}
\newtheorem{remark}[theorem]{Remark}

\newcommand{\supp}{\operatorname{supp}}
\newcommand{\R}{\mathbb R}

\DeclareMathOperator{\dom}{dom}

\title{Weighted fixed-bidegree subelliptic estimates for the \(\overline{\partial}_f\)-Neumann problem on Lipschitz pseudoconvex domain}
\author{Ping'an Zhang\footnote{\noindent{School of Mathematics and Statistics, Xidian University, 266 Xinglong Section of Xifeng Road, Xi'an, Shaanxi, 710126, P. R. China.}
\noindent{E-mail addresses: pinganzhang@xidian.edu.cn}}}
\date{}

\begin{document}
\maketitle
\begin{abstract}
We establish two types of \emph{a priori} subelliptic estimates for the twisted \(\overline{\partial}_f\)-Neumann problem on a bounded pseudoconvex domain with Lipschitz boundary. The estimates feature a subelliptic gain in \((0, 1/2)\) and are uniform with respect to the bidegree. The boundary is assumed to satisfy a plurisubharmonic-type condition, characterized by
the local existence, near each boundary point, of a bounded plurisubharmonic function whose Hessian blows up algebraically near the boundary, with potentially different exponents at different points.
The twisting function is holomorphic in a neighborhood of the domain, has bounded gradient, and its critical set inside the domain consists of isolated interior points with finite vanishing orders.
\end{abstract}



\section{Introduction and Global Assumptions}
\label{sec:assumptions}

A subelliptic estimate for the \(\overline{\partial}\)-Neumann problem asserts that the solution has a fractional gain of derivatives. On smoothly bounded strictly pseudoconvex domains, the classical basic estimate yields such an estimate with \(\epsilon=1/2\); see \cite{ChenShaw,Kohn1963,Kohn1964,Kohn1979}, and \cite[Corollary~3.2.12 and Theorem~2.4.4]{FollandKohn1972} and \cite[Proposition~3.1]{Straube2010}. 

In the smooth pseudoconvex case, subellipticity at a boundary point is
equivalent to finite type in the appropriate \(q\)-dimensional sense:
for \((p,q)\)-forms ($1\le q<n$), a subelliptic estimate of some order
\(\epsilon>0\) holds near a boundary point \(z_0\) if and only if Catlin's
\(q\)-type \(D_q(z_0)\) is finite, i.e.\ the generic order of contact of \(q\)-dimensional complex varieties with the boundary at \(z_0\) is bounded \cite[Theorem~1.1]{Catlin1987}. The sufficiency of the geometric condition \(D_q(z_0)<\infty\) is proved in \cite{Catlin1987}, while necessity is due to Catlin's earlier work (see \cite[Theorem~3.3]{Catlin1987}).  When \(q=1\), \(D_1\) coincides with the
ordinary D'Angelo \(1\)-type; for \(q>1\) the relevant invariant is Catlin's
\(q\)-type \(D_q(z_0)\), not the ordinary \(1\)-type (its agreement with
D'Angelo's \(q\)-type is left open in \cite[Section~1 and the remark after Theorem~3.3]{Catlin1987}).
This characterization is formulated for smooth boundaries, and the boundary regularity is essential to the classical formulation.
Catlin's weighted estimate \cite[Theorem~2.2]{Catlin1987} shows, for a smoothly bounded domain that is pseudoconvex near \(z_0\), that if, for every small \(\delta>0\), there is a \(C^2\) function \(\lambda_\delta\) on a neighborhood \(V\) of \(z_0\) with \(|\lambda_\delta|\le1\) on \(V\), whose complex Hessian, acting on \((p,q)\)-forms, is nonnegative on \(V\) and bounded below by \(c\,\delta^{-2\epsilon}\) (with \(c\) independent of \(\delta\)) on \(V\cap\{-\delta<r\le0\}\), then an \(\epsilon\)-subelliptic estimate holds on some neighborhood \(V'\Subset V\) of \(z_0\).

Derridj characterized the maximal estimate for \((0,1)\)-forms by the Levi eigenvalues: for a neighborhood \(V\) of a boundary point, the maximal estimate holds on \(V\) if and only if there is \(A>0\) such that every eigenvalue of the Levi form is bounded from below by \(A\) times the trace of the Levi form on \(\partial\Omega\cap V\) \cite[Th\'eor\`eme~3.1]{Derridj1978}. He then proved that a bounded pseudoconvex domain with real-analytic boundary satisfying this condition near \(\partial\Omega\) satisfies a subelliptic estimate for \((0,1)\)-forms at every boundary point \cite[Th\'eor\`eme~4.1]{Derridj1978}; in particular, every bounded pseudoconvex domain in \(\mathbb C^2\) with real-analytic boundary is subelliptic \cite[Corollaire~4.4]{Derridj1978}. Ben Moussa later proved a subelliptic estimate for \((0,q)\)-forms on a neighborhood of a connected component \(D\) of the Levi-degeneracy set under part of the hypotheses of his own semi-global analyticity theorem (real-analyticity of the boundary near \(D\), a tangential frame, and a Derridj-type lower bound for sums of \(q\) Levi eigenvalues by the trace) \cite[Proposition~3.9]{BenMoussa2000}, and in another setting assumed subellipticity directly, together with a real-analytic frame and an exact commutation hypothesis, to obtain semi-global analytic hypoellipticity \cite[Th\'eor\`eme~3.4]{BenMoussa2000}. 

For piecewise smooth strictly pseudoconvex boundaries, Michel and Shaw \cite[Theorem~1.2, Corollaries~1.3--1.4]{MichShaw} establish \(H^{1/2}\) estimates for \(N_q\), \(\overline{\partial} N_q\), and \(\overline{\partial}^{\ast} N_q\) (\(1\le q\le n-1\)), compactness of these operators, and the corresponding a priori estimate, by exploiting the boundary term in the Kohn--Morrey inequality on smooth approximating domains; they indicate that the same estimates hold for strictly pseudoconvex domains with Lipschitz boundary \cite[Theorem~3.1]{MichShaw}. For bounded pseudoconvex domains, Straube \cite{Straube1997} analyzed a purely pluripotential-theoretic sufficient condition: if near a boundary point \(P\) the boundary is the graph of a Lipschitz function and there is one bounded plurisubharmonic function on \(U\cap\Omega\) (\(U\) a neighborhood of \(P\)) whose Levi form dominates \(\operatorname{const.}\,\delta^{-2\epsilon}\) times the Euclidean Levi form, with \(0\le\epsilon<1/2\), then an \(\epsilon\)-subelliptic estimate holds near \(P\) \cite[Theorem~1]{Straube1997}.
In particular, the class of piecewise smooth boundaries of finite type covered by Straube's theorem \cite[Theorem~2]{Straube1997}---namely, local transverse intersections of finitely many smooth finite-type hypersurfaces, each pseudoconvex from the indicated side---satisfies this hypothesis, so the \(\overline{\partial}\)-Neumann problem is subelliptic near such boundary points.

The introduction of the twisted differential operator \(\overline\partial_f := \overline\partial + \partial f \wedge\), along with its formal adjoint \(\vartheta_f := \vartheta + \iota_{(\partial f)^\sharp}\) and its Hilbert-space adjoint \(\overline\partial_f^{\ast}\) (defined via the graded-adjoint convention \((\overline\partial_f^{\ast}u,v)=(u,\overline\partial_f v)\)), is of foundational importance in the rigorous deformation theory of supersymmetric quantum mechanics. The significance of these operators arises directly from the holomorphic function \(f\), which serves as the superpotential. In the physical framework of the \(\mathcal N=2\) Landau-Ginzburg model, the superpotential \(f\) completely determines the chiral ring and dictates the pseudo-topological ground state metric (the \(tt^{\ast}\) metric, which in the conformal limit essentially reduces to the Zamolodchikov metric) obtained by fusing the topological theory with its anti-topological conjugate, see Cecotti and Vafa~\cite{CecoVafa}. For $\mathcal N=2$ supersymmetric models, see Cecotti \cite{Cecotti}. Geometrically and analytically, the introduction of \(f\) fundamentally twists the standard Dolbeault complex on the K\"{a}hler manifold \((M,g)\), thereby generating the twisted Laplace operator \(\Box_f =  \overline\partial_f \overline\partial_f^\ast  + \overline\partial_f^\ast \overline\partial_f\).  Consequently, studying the deformations of these twisted operators with respect to \(f\) yields the Cecotti-Vafa equations (the \(tt^{\ast}\) equations). These equations encode the integrable structure, flat connections, and generalized special geometry that are essential for understanding the deformation space of the Landau-Ginzburg B-model \cite{CecoVafa, Fan}. The motivation for the present work stems from the study in Wen and Fan \cite{WeFan}. They make a sharper estimation for studying the global regularity for the $\overline\partial_f$-Neumann operator.

Our main focus is the effect of holomorphic twisted functions with isolated interior critical points. In both the unweighted and weighted regimes, the boundary contributions are handled by invoking the known subelliptic estimates of Straube \cite{Straube1997}. For the unweighted estimate (Theorem~\ref{them:global}), combining this boundary control with the interior \(H^1\) estimate via the bidegree bridge \eqref{eq:bridge} shows that interior critical points cause no loss, since the twist is zero-order; in this estimate \(f\) enters only through the boundedness of \(\partial f\) on \(\Omega\). In the weighted Sobolev scale (Theorem~\ref{thm:global-Hs}), the weight \(\theta\) vanishes to order \(m_a\) at each critical point \(z_a\), so that the commutator \(\overline\partial\theta\wedge u\) is not dominated by \(\theta u\) there; this loss is compensated by the Hardy-type inequalities of Section~\ref{sec:critical}, and we prove that the twisted \(\overline{\partial}_f\)-Neumann system still satisfies an \(\epsilon\)-subelliptic estimate.

Throughout, \(\mathbb{C}^n\) is equipped with its standard Hermitian metric and the associated flat Chern connection, and we assume the following.

\begin{enumerate}[label=\textup{(\roman*)}]

\item Lipschitz pseudoconvex domain.
Let \(\Omega\Subset\mathbb C^n\) be a bounded domain with boundary
\(\partial\Omega\), and let \(U_\rho\) be an open neighborhood of
\(\overline\Omega\) carrying a real-valued Lipschitz defining function
\(\rho:U_\rho\to\mathbb R\) such that
\[
\Omega=\{z\in U_\rho:\rho(z)<0\},\qquad
\partial\Omega=\{z\in U_\rho:\rho(z)=0\}.
\]
Assume that \(\Omega\) is pseudoconvex, i.e.\ \(-\log\operatorname{dist}(\cdot,\partial\Omega)\)
is plurisubharmonic on \(\Omega\).

\item {Boundary geometry hypothesis.}
 As in Straube \cite{Straube1997}, we write
\(d^c=i(\overline\partial-\partial)\), so that \(dd^c=2i\partial\overline\partial\).
For every boundary point \(P\in\partial\Omega\), assume that the boundary
is the graph of a Lipschitz function near \(P\) and that there exist a
neighborhood \(U_P\) of \(P\), a number \(\epsilon_P\in(0,1/2)\), and a
bounded plurisubharmonic function \(\lambda_P\) on \(U_P\cap\Omega\) such
that
\begin{equation}\label{eq:boundary-hyp}
dd^c\lambda_P\ge \operatorname{const.}\delta^{-2\epsilon_P}dd^c|z|^2
\quad\text{in }U_P\cap\Omega
\end{equation}
as currents.  The exponent \(\epsilon_P\) is constant on the whole
neighborhood \(U_P\).  Here \(\delta(z)=\operatorname{dist}(z,\partial\Omega)\)
and the constant \(\operatorname{const.}>0\) is independent of \(z\). This condition is adopted from Ref. \cite{Straube1997}.

Since \(P\in\partial\Omega\), by taking \(U_P\) sufficiently small we may
assume in addition that
\[
\delta(z)\le 1
\qquad\text{for all }z\in U_P\cap\Omega .
\]
Indeed, since \(P\in\partial\Omega\), we have \(\delta(z)\le|z-P|\); it thus suffices to take \(U_P\subset B_1(P)\). Shrinking \(U_P\) only restricts the set on which the inequality \eqref{eq:boundary-hyp} is required, so the lower bound remains valid on the smaller neighborhood.

\item {Holomorphic twisting function.}
After possibly shrinking \(U_\rho\) such that it still contains \(\overline\Omega\), assume that \(f\in\mathcal O(U_\rho)\).  Put
$s:=\partial f|_{\Omega}$,
a holomorphic \((1,0)\)-form on \(\Omega\).  Since \(\overline{\Omega}\)
is compact and \(f\) is holomorphic on \(U_\rho\), we have
\[
\partial f\in C^\infty(\overline{\Omega};\Lambda^{1,0}),
\qquad
\|\partial f\|_{L^\infty(\Omega;\Lambda^{1,0})}<\infty .
\]

\item {Isolated critical set and finite vanishing orders.}
\label{ass:order}
Define the critical set
$Z:=\{z\in\Omega:s(z)=0\}$.
We assume that every point of \(Z\) is an isolated zero of \(s\) and that
\(Z\neq\varnothing\).  Then \(Z\) is finite.  Indeed,
\(A:=\{z\in U_\rho:\partial f(z)=0\}\) is an analytic subset of \(U_\rho\); a point
of \(Z\) that is isolated in \(Z\) lies on no positive-dimensional
irreducible component of \(A\) (such a component meets the open set
\(\Omega\) in a set without isolated points), so it is a zero-dimensional
component of \(A\); since the irreducible components of \(A\) are locally
finite in \(U_\rho\) and \(\overline\Omega\subset U_\rho\) is compact, only finitely
many of them meet \(\overline\Omega\).  In particular \(Z\) is a finite
subset of the open set \(\Omega\), hence \(Z\Subset\Omega\).  We write \(z_c\) for a
generic point of \(Z\), and use the index letter \(a\) as the summation index
over \(Z\), with \(z_a\in Z\), \(a=1,\dots,N\), $N\ge 1$.
For each \(z_c\in Z\), the vanishing order
$m_c:=\operatorname{ord}_{z_c}(s)=\operatorname{ord}_{z_c}(\partial f)$
is defined componentwise as the minimum of the vanishing orders of the
coefficients of \(s\) in a local frame; it is finite and satisfies
\(m_c\ge1\). The two conditions (iii) and (iv) have been studied in \cite{WeFan}.

\end{enumerate}

Set the distance function
$\upsilon(z):=d(z,Z)$.
Since \(Z\Subset\Omega\), the distance
\(d(Z,\partial\Omega)>0\).

\emph{Order-adapted weight.}
Choose \(R_a>0\) for each \(z_a\in Z\) so small that
$B_{2R_a}(z_a)\Subset\Omega$,
the balls \(B_{2R_a}(z_a)\) are pairwise disjoint, and \(z_a\) is the unique
nearest point of \(Z\) to every \(z\in B_{2R_a}(z_a)\); equivalently,
\[
\upsilon(z)=|z-z_a|=:r_a(z)
\qquad
\text{for } z\in B_{2R_a}(z_a).
\]
Such \(R_a\) exist because \(Z\) is finite and \(Z\Subset\Omega\). For
instance, if \(N\ge2\), it suffices to require, for every \(a\),
\[
4R_a<\min_{b\ne a}|z_a-z_b|
\quad\text{and}\quad
2R_a<d(Z,\partial\Omega).
\]
If \(N=1\), the pairwise-disjointness condition is vacuous and the unique-nearest-point condition is automatically satisfied; in this case it suffices to require
$2R_1<d(Z,\partial\Omega)$.

The weight exponent used near \(z_a\) is the order \(m_a\) itself; we
call the resulting weight the \emph{order-adapted weight}.
Fix once and for all a smooth cutoff
\(\zeta:[0,\infty)\to[0,1]\) with \(\zeta\equiv1\) on \([0,1]\) and
\(\zeta\equiv0\) on \([2,\infty)\). Define
\(\theta:\overline\Omega\to[0,\infty)\) by

\[
\theta(z)=
\begin{cases}
r_a(z)^{m_a}, & z\in B_{R_a}(z_a),\\[1mm]
r_a(z)^{m_a}\zeta\!\left(\dfrac{r_a(z)}{R_a}\right)
+r_a(z)\left(1-\zeta\!\left(\dfrac{r_a(z)}{R_a}\right)\right),
& z\in B_{2R_a}(z_a)\setminus B_{R_a}(z_a),\\[1mm]
\upsilon(z), & z\in \overline\Omega\setminus\textstyle\bigcup_aB_{2R_a}(z_a).
\end{cases}
\]
At \(r_a=R_a\) we have \(\zeta=1\), so the transition formula equals
\(r_a^{m_a}\); at \(r_a=2R_a\) we have \(\zeta=0\), so it equals \(r_a\).
Since \(\upsilon=r_a\) on \(B_{2R_a}(z_a)\), the value on the gluing sphere
\(r_a=2R_a\) is \(2R_a\), which matches the exterior value
\(\upsilon=2R_a\).
Hence \(\theta\) is continuous and positive on
\(\overline\Omega\setminus Z\). 
On the transition annulus
\(R_a\le r_a\le2R_a\) 
the radius is bounded away from \(0\) and \(\infty\),
so \(\theta\) is a smooth positive function bounded above and below by
positive constants depending only on \(R_a\) and \(m_a\); its derivative is
bounded there as well.
Outside \(\bigcup_aB_{2R_a}(z_a)\), \(\theta=\upsilon\)
is the distance to the finite set \(Z\). 

Because \(Z\Subset\Omega\), we have
\[
d_0:=\min_a 2R_a>0,\qquad d_0\le\upsilon(z)\le\operatorname{diam}(\Omega)
\quad\text{for } z\in\overline\Omega\setminus\textstyle\bigcup_aB_{2R_a}(z_a),
\]
and \(|\nabla\upsilon|\le1\) almost everywhere. Consequently \(\theta\) is
Lipschitz on \(\overline\Omega\) (it is Lipschitz on each of the three
pieces and continuous across the gluing spheres), so
\(\theta\in W^{1,\infty}(\Omega)\);
and on every boundary collar and on
\(\overline\Omega\setminus\bigcup_aB_{2R_a}(z_a)\) the weight \(\theta\) is bounded
above and below by positive constants.
The constant \(d_0\) depends on the
choice of \(R_a\); it is fixed once the balls are chosen, and set \(R_a=R_c\) when \(z_a=z_c\).

\emph{Full exterior algebra and \(L^2\)-spaces.}
For \(0\le k\le 2n\), set
\[
\Lambda^k:=\bigoplus_{p+q=k}\Lambda^{p,q},
\qquad
L^2_{\Lambda^k}(\Omega)
:=
\bigoplus_{p+q=k}L^2_{(p,q)}(\Omega),
\qquad
C^{\infty}_{\Lambda^k}(\Omega)
:=
\bigoplus_{p+q=k}C^{\infty}_{(p,q)}(\Omega).
\]
We extend the grading by the endpoint conventions
\[
\Lambda^{-1}:=\Lambda^{2n+1}:=\{0\},
\qquad
L^2_{\Lambda^{-1}}(\Omega):=L^2_{\Lambda^{2n+1}}(\Omega):=\{0\}.
\]
When smooth endpoint spaces are needed, we also set
\[
C^\infty_{\Lambda^{-1}}(\Omega)
:=C^\infty_{\Lambda^{2n+1}}(\Omega):=\{0\}.
\]

The \(L^2\)-inner product on \(L^2_{\Lambda^k}(\Omega)\) is the direct sum of
the Hermitian \(L^2\)-inner products on the summands \(L^2_{(p,q)}(\Omega)\).
We also write
\[
L^2_{\Lambda^\bullet}(\Omega)
:=
\bigoplus_{k=0}^{2n}L^2_{\Lambda^k}(\Omega).
\]

\emph{Inner-product convention.}
We use the Hermitian inner product on \(\Lambda^{p,q}\mathbb C^n\) for
which the constant coframe
\(\{dz_I\wedge d\overline{z}_J:|I|=p,\ |J|=q\}\) is
orthonormal, i.e.\ \(|dz_i|=|d\overline{z}_i|=1\), and \(L^2\)-norms are taken
with respect to Lebesgue measure \(dV\) on \(\mathbb C^n=\mathbb R^{2n}\).
Write \(\Delta:=\sum_{j=1}^n(\partial_{x_j}^2+\partial_{y_j}^2)\) for the
real Euclidean Laplacian and \(\nabla\) for the real Euclidean gradient,
acting coefficientwise in the constant coframe.  With this convention
\(\partial_{z_j}\partial_{\overline{z}_j}=\frac14(\partial_{x_j}^2+\partial_{y_j}^2)\),
and the flat identity reads
\[
(\overline\partial+\overline\partial^{\ast})^2
=\overline\partial\,\overline\partial^{\ast}+\overline\partial^{\ast}\overline\partial
=-\tfrac14\Delta
\qquad\text{coefficientwise on }\Lambda^{\bullet,\bullet}\mathbb C^n .
\]
Consequently, for compactly supported smooth forms \(\varphi\) (of any,
possibly mixed, bidegree),
\[
\|\nabla\varphi\|_{L^2}^2
=4\|(\overline\partial+\overline\partial^{\ast})\varphi\|_{L^2}^2
=4\|\overline\partial\varphi\|_{L^2}^2+4\|\overline\partial^{\ast}\varphi\|_{L^2}^2 ,
\]
where the second equality uses
\((\overline\partial\varphi,\overline\partial^{\ast}\varphi)=(\overline\partial^2\varphi,\varphi)=0\).
Explicitly, with \(e_j:=d\overline z_j\wedge\) and its pointwise adjoint
\(e_j^{\ast}=\iota_{(d\overline z_j)^\sharp}\),
\[
\mathcal D:=\overline\partial+\overline\partial^{\ast}
=\sum_{j=1}^n\bigl(e_j\partial_{\overline z_j}-e_j^{\ast}\partial_{z_j}\bigr)
=\frac12\sum_{j=1}^n\bigl(c_{x_j}\partial_{x_j}+c_{y_j}\partial_{y_j}\bigr),
\qquad
c_{x_j}:=e_j-e_j^{\ast},\quad c_{y_j}:=i(e_j+e_j^{\ast}).
\]
In this paragraph and in Lemma~\ref{lem:angular-spectrum-gap}, the
letters \(a,b,d\) used as subscripts of \(c_a\), \(x_a\), \(X_{ab}\) run
over the real coordinate indices \(\{1,\dots,2n\}\); they are unrelated
to the index \(a\) labelling the points \(z_a\in Z\).
From \(e_je_k^{\ast}+e_k^{\ast}e_j=\delta_{jk}\) and
\(e_je_k+e_ke_j=0\) one checks that the \(2n\) constant matrices
\(c_a\in\{c_{x_j},c_{y_j}\}\) are skew-adjoint and satisfy
\(c_ac_b+c_bc_a=-2\delta_{ab}\); in particular each \(c_a\) is unitary.
Thus \(\mathcal D=\frac12\mathfrak D\), where
\(\mathfrak D=\sum_{a=1}^{2n}c_a\partial_{x_a}\) is a Euclidean Dirac
operator on the Clifford module \(\Lambda^{\bullet,\bullet}\mathbb C^n\),
and the principal symbol of \(\mathcal D\) at a real covector \(\xi\) has
norm \(|\xi|/2\).
All constants below are computed in this convention.  

We denote by \(\|\cdot\|\) the standard \(L^2(\Omega)\) norm (or the \(L^2(\mathbb{C}^n)\) norm for forms extended by zero outside \(\Omega\)), unless a subscript is explicitly specified.

\emph{Twisted operators as full exterior-algebra operators.}
Let \(\vartheta\) denote the formal adjoint of \(\overline\partial\), so that
\(\vartheta\) is a constant-coefficient differential operator characterized by
\[
(\overline\partial u,v)=(u,\vartheta v)
\]
for all compactly supported smooth forms \(u,v\) of appropriate degrees.
Let \(\overline\partial_f\) be the maximal \(L^2\)-realization of
\(\overline\partial+\partial f\wedge\), and let
$\vartheta_f:=\vartheta+\iota_{(\partial f)^\sharp}$ denote the formal adjoint of
\(\overline\partial+\partial f\wedge\).
Let
\(\overline\partial_f^{\ast}\) denote its Hilbert-space adjoint in the
graded-adjoint convention:
\[
(\overline\partial_f^{\ast}u,v)=(u,\overline\partial_f v)
\]
for all \(u\in\dom(\overline\partial_f^{\ast})\) and
\(v\in\dom(\overline\partial_f)\).
For compactly supported smooth forms \(v\in C^\infty_{c,\Lambda^k}(\Omega)\),
which belong to \(\dom(\overline\partial_f^{\ast})\), the Hilbert adjoint is
given by this formal expression:
\[
\overline\partial_f^{\ast}v=\vartheta_f v
=\vartheta v+\iota_{(\partial f)^\sharp}v.
\]
Consequently,
\[
(\overline\partial_f u,v)=(u,\overline\partial_f^{\ast}v)
=(u,\vartheta_f v)
\]
for all compactly supported smooth forms \(u,v\) of appropriate degrees.
Thus, on compactly supported smooth forms, the formal adjoint
\(\vartheta_f\) and the Hilbert adjoint \(\overline\partial_f^{\ast}\)
coincide.

Here \(\partial f\) is a \((1,0)\)-form and \((\partial f)^\sharp\) is the
associated \((1,0)\)-vector field. We use the endpoint convention \(\Lambda^{-1}:=\{0\}\). Then both summands of \(\vartheta_f=\vartheta+\iota_{(\partial f)^\sharp}\) vanish on \(\Lambda^0\), so that \(\vartheta_f=0\) on \(\Lambda^0\).
On a fixed bidegree summand \(\Lambda^{p,q}\), the two summands of
\(\overline{\partial}_fu\) lie in the distinct bidegrees
\[
\overline\partial u\in\Lambda^{p,q+1},
\qquad
\partial f\wedge u\in\Lambda^{p+1,q},
\]
and the two summands of \(\overline{\partial}_f^{\ast}u\) lie in the
distinct bidegrees
\[
\overline\partial^{\ast} u\in\Lambda^{p,q-1},
\qquad
\iota_{(\partial f)^\sharp}u\in\Lambda^{p-1,q}.
\]
This bidegree separation is asserted only for a fixed-bidegree form
\(u\in\Lambda^{p,q}\).
For a mixed-bidegree form
\[
u=\sum_{p+q=k}u_{p,q}\in\Lambda^k,
\]
the contributions from different bidegree components may occupy the same
bidegree; in particular, no pointwise orthogonality or Pythagorean
identity for \(\overline{\partial}_fu\) or
\(\overline{\partial}_f^{\ast}u\) is claimed at the level of the full
exterior algebra.  Such identities are used below only after restricting
to a fixed bidegree \((p,q)\); for mixed-bidegree objects only the Dirac
energy and estimates that do not require a bidegree comparison are used.

\emph{\(L^2\)-realizations and the weak Neumann condition.}
Let
\[
\overline{\partial}_f:=\bigoplus_{k=0}^{2n}\overline{\partial}_{f,k},
\qquad
\overline{\partial}_{f,k}:\dom(\overline{\partial}_{f,k})\subset
L^2_{\Lambda^k}(\Omega)\to L^2_{\Lambda^{k+1}}(\Omega),
\]
where, for \(0\le k<2n\), \(\overline{\partial}_{f,k}\) is the maximal
\(L^2\)-realization of \(\overline\partial+\partial f\wedge\) on
\(L^2_{\Lambda^k}\), i.e.\ its domain consists of those
\(u\in L^2_{\Lambda^k}(\Omega)\) for which the distribution
\(\overline\partial u+\partial f\wedge u\) lies in
\(L^2_{\Lambda^{k+1}}(\Omega)\). At the top endpoint we set
\[
\overline{\partial}_{f,2n}:=0,
\qquad
\dom(\overline{\partial}_{f,2n}):=L^2_{\Lambda^{2n}}(\Omega).
\]
This agrees with the maximal-realization definition because
\(L^2_{\Lambda^{2n+1}}(\Omega)=\{0\}\), and for
\(u\in\Lambda^{2n}\) both \(\overline\partial u\) and
\(\partial f\wedge u\) vanish.

For \(u\in L^2_{\Lambda^k}(\Omega)\) we write
\(\overline{\partial}_fu:=\overline{\partial}_{f,k}u\).
We use the graded-adjoint convention:
\[
\overline{\partial}_f^{\ast}:=\bigoplus_{k=0}^{2n}\overline{\partial}_{f,k}^{\ast},
\qquad
\overline{\partial}_{f,k}^{\ast}:\dom(\overline{\partial}_{f,k}^{\ast})\subset
L^2_{\Lambda^k}(\Omega)\to L^2_{\Lambda^{k-1}}(\Omega),
\]
where \(\overline{\partial}_{f,k}^{\ast}\) is the Hilbert-space adjoint of
\(\overline{\partial}_{f,k-1}\). Here
\[
\overline{\partial}_{f,-1}:=0:
L^2_{\Lambda^{-1}}(\Omega)\to L^2_{\Lambda^0}(\Omega),
\]
so that for \(k=0\),
\[
\overline{\partial}_{f,0}^{\ast}=0,
\qquad
\dom(\overline{\partial}_{f,0}^{\ast})
=L^2_{\Lambda^0}(\Omega).
\]
\(u\in\dom(\overline{\partial}_{f,k}^{\ast})\) if and only if there
is \(w\in L^2_{\Lambda^{k-1}}(\Omega)\) with
\((u,\overline{\partial}_fv)=(w,v)\) for all \(v\in\dom(\overline{\partial}_{f,k-1})\), and then
\(\overline{\partial}_{f,k}^{\ast}u=w\).  The untwisted operators
\(\overline\partial\), \(\overline\partial^{\ast}\) on \(L^2_{\Lambda^\bullet}(\Omega)\) are
defined in the same way with \(f=0\). In particular, with the same endpoint conventions,
\[
\overline\partial_{2n}=0,
\qquad
\overline\partial_0^{\ast}=0.
\] 
On a single bidegree, for \(0\le p\le n\), \(0\le q\le n\),
\[
\overline\partial_{(p,q)}:L^2_{(p,q)}(\Omega)\to L^2_{(p,q+1)}(\Omega)
\]
denotes the maximal realization, with the endpoint conventions
\[
L^2_{(p,-1)}(\Omega):=L^2_{(p,n+1)}(\Omega):=\{0\},
\]
\[
\overline\partial_{(p,n)}:=0,
\qquad
\dom(\overline\partial_{(p,n)}):=L^2_{(p,n)}(\Omega),
\]
and
\[
\overline\partial_{(p,-1)}:=0.
\]
The operator \(\overline\partial^{\ast}_{(p,q)}\) is the Hilbert
adjoint of \(\overline\partial_{(p,q-1)}\). In particular,
\[
\overline\partial^{\ast}_{(p,0)}=0,
\qquad
\dom(\overline\partial^{\ast}_{(p,0)})
=L^2_{(p,0)}(\Omega).
\]
Since \(\partial f\wedge\partial f=0\) and \(\overline\partial\partial f=0\),
one has \(\overline{\partial}_f\circ\overline{\partial}_f=0\) in the sense of distributions.

\begin{lemma}[Domains of the twisted operators and localization]
\label{lem:domains}
\begin{enumerate}[label=\textup{(\alph*)}]
\item For every \(k\),
\(\dom(\overline{\partial}_{f,k})=\dom(\overline\partial_k)\) and
\(\dom(\overline{\partial}_{f,k}^{\ast})=\dom(\overline\partial_k^{\ast})\), with
\(\overline{\partial}_fu=\overline\partial u+\partial f\wedge u\) and
\(\overline{\partial}_f^{\ast}u=\overline\partial^{\ast}u+\iota_{(\partial f)^\sharp}u\).
\item Let \(u=\sum_{p+q=k}u_{p,q}\in L^2_{\Lambda^k}(\Omega)\) with
\(u_{p,q}\in L^2_{(p,q)}(\Omega)\).  Then \(u\in\dom(\overline\partial_k)\) if
and only if \(u_{p,q}\in\dom(\overline\partial_{(p,q)})\) for all \((p,q)\),
and \(u\in\dom(\overline\partial_k^{\ast})\) if and only if
\(u_{p,q}\in\dom(\overline\partial^{\ast}_{(p,q)})\) for all \((p,q)\); in that
case \(\overline\partial u=\sum\overline\partial_{(p,q)}u_{p,q}\) and
\(\overline\partial^{\ast}u=\sum\overline\partial^{\ast}_{(p,q)}u_{p,q}\).
\item Let \(u\in L^2_{\Lambda^k}(\Omega)\) with \(\supp u\) compact in
\(\Omega\), and let \(\vartheta\) denote the formal adjoint of
\(\overline\partial\) (a constant-coefficient operator).  Then
\(u\in\dom(\overline\partial^{\ast})\) if and only if the distribution
\(\vartheta u\) belongs to \(L^2\), and then \(\overline\partial^{\ast}u=\vartheta u\).
In particular, for such \(u\) and any open \(G\) with
\(\supp u\subset G\subset\Omega\), membership of \(u\) in
\(\dom(\overline\partial)\cap\dom(\overline\partial^{\ast})\) (equivalently, by
(a), in \(\dom(\overline{\partial}_f)\cap\dom(\overline{\partial}_f^{\ast})\)) is the same whether these
domains are formed on \(\Omega\) or on \(G\).
\end{enumerate}
\end{lemma}

\begin{proof}
(a) The coefficients of \(\partial f\) are bounded on \(\overline\Omega\),
so \(\mathcal B:=\partial f\wedge\) is a bounded, everywhere defined operator on
\(L^2_{\Lambda^\bullet}(\Omega)\) whose Hilbert adjoint is the pointwise
adjoint \(\iota_{(\partial f)^\sharp}\).  The distribution
\(\overline\partial u+\mathcal Bu\) lies in \(L^2\) if and only if \(\overline\partial u\) does,
which gives the first identity.  For a densely defined operator \(S\) and
a bounded everywhere defined \(\mathcal B\), \((S+\mathcal B)^{\ast}=S^{\ast}+\mathcal B^{\ast}\) with
\(\dom((S+\mathcal B)^{\ast})=\dom(S^{\ast})\), since \(v\mapsto(u,\mathcal Bv)\) is bounded
for every \(u\).  Applying this with \(S=\overline\partial_{k-1}\) gives the second
identity, see also \cite[page 67]{WeFan}.

(b) The operator \(\overline\partial_{k-1}\) is the finite orthogonal direct sum
\(\bigoplus_{p+q=k-1}\overline\partial_{(p,q)}\) acting from
\(\bigoplus L^2_{(p,q)}\) to \(\bigoplus L^2_{(p,q+1)}\).  The statement for
\(\dom(\overline\partial_k)\) is immediate since \(\overline\partial\) preserves the
holomorphic degree.  For a finite orthogonal direct sum of densely defined
operators \(A_j\), \((\bigoplus_jA_j)^{\ast}=\bigoplus_jA_j^{\ast}\) with
domain \(\bigoplus_j\dom(A_j^{\ast})\), as one sees by testing the
adjoint relation with families having a single nonzero entry. 
For every \(0\le p\le n\), the bidegree-\((p,0)\) component is not in the range of any \(\overline\partial_{(p,q)}\), since \(\overline\partial_{(p,q)}\) raises the antiholomorphic degree; this accounts for the convention \(\overline\partial^{\ast}_{(p,0)}=0\). In the total-degree notation, the \((k,0)\) component occurs only when \(k\le n\).
At the top antiholomorphic degree, the convention
\(\overline\partial_{(p,n)}=0\) is handled in the same way; it does not
affect the finite direct-sum decomposition because
\(L^2_{(p,n+1)}(\Omega)=\{0\}\).

(c) If \(u\in\dom(\overline\partial^{\ast})\), testing with
\(v\in C^\infty_c\) gives \(\vartheta u=\overline\partial^{\ast}u\in L^2\).

Conversely, assume \(\vartheta u\in L^2\). 
Choose a cutoff function \(\chi\in C^\infty_c(\Omega)\)
with \(\chi\equiv1\) near \(\supp u\).
For \(v\in\dom(\overline\partial)\),
\(\chi v\in\dom(\overline\partial)\) with
\(\overline\partial(\chi v)=\chi\overline\partial v+\overline\partial\chi\wedge v\)
(Leibniz rule for a smooth multiplier), and \(\overline\partial\chi=0\) near
\(\supp u\); hence \((u,\overline\partial v)=(u,\overline\partial(\chi v))\).  The mollifications
\(J_\eta\ast(\chi v)\in C^\infty_c(\Omega)\) converge to \(\chi v\) in
the graph norm of \(\overline\partial\), because \(\overline\partial\) has constant
coefficients.  Therefore
\((u,\overline\partial(\chi v))=\lim_\eta(u,\overline\partial(J_\eta\ast(\chi v)))
=\lim_\eta(\vartheta u,J_\eta\ast(\chi v))=(\vartheta u,\chi v)=(\vartheta u,v)\),
the last step because \(\chi\equiv1\) on \(\supp u\supset\supp\vartheta u\).
Thus \(u\in\dom(\overline\partial^{\ast})\) with \(\overline\partial^{\ast}u=\vartheta u\).  The
same argument on \(G\) in place of \(\Omega\) gives the final assertion,
since both conditions reduce to \(\overline\partial u,\vartheta u\in L^2\).
\end{proof}

We refer to \(u\in\dom(\overline{\partial}_f^{\ast})\) as the \emph{weak \(\overline{\partial}_f\)-Neumann
condition}.  By Lemma~\ref{lem:domains}(a),(b) it is, for every bidegree,
exactly the weak \(\overline\partial\)-Neumann condition
\(u_{p,q}\in\dom(\overline\partial^{\ast}_{(p,q)})\), which is the variational
meaning of \(\iota_{(\overline{\partial}\rho)^\sharp}u=0\) on \(\partial\Omega\) in the Lipschitz setting.
In particular
\begin{equation}
\label{eq:domain-equality}
\dom(\overline{\partial}_f)=\dom(\overline\partial),
\qquad
\dom(\overline{\partial}_f^{\ast})=\dom(\overline\partial^{\ast})
\qquad\text{on every }\Lambda^{p,q},
\end{equation}
and the comparison between the twisted and untwisted operators used
below rests on \eqref{eq:domain-equality}, on
Lemma~\ref{lem:domains}(b),(c), and on the bridge \eqref{eq:bridge}
(established after Lemma~\ref{lem:componentwise}).  When \(u\in L^2_{(p,q)}(\Omega)\), the notation
\(u\in\dom(\overline{\partial}_f)\cap\dom(\overline{\partial}_f^{\ast})\) always refers to these total-degree
realizations; in particular \(\overline{\partial}_fu\in\Lambda^{p,q+1}\oplus\Lambda^{p+1,q}\)
and \(\overline{\partial}_f^{\ast}u\in\Lambda^{p,q-1}\oplus\Lambda^{p-1,q}\).

The following result is abstracted from the argument in \cite[Proposition~2.8]{Straube2010}.

\begin{lemma}[Abstract Hilbert-complex Laplacian domains]
\label{lem:abstract-Hilbert-complex}
Let \(\mathrm{H}^\bullet=\bigoplus_k \mathrm{H}^k\) be a graded Hilbert space.  For each grading index
\(k\), let
\[
T_k:\dom(T_k)\subset \mathrm{H}^k\longrightarrow \mathrm{H}^{k+1}
\]
be a closed, densely defined operator, and let
\[
T_k^{\ast}:\dom(T_k^{\ast})\subset \mathrm{H}^{k+1}\longrightarrow \mathrm{H}^k
\]
be its Hilbert adjoint.  We assume that
$T_{k+1}T_k=0$
in the strong sense that
\[
\dom(T_k)\subset\dom(T_{k+1}T_k),
\qquad
T_{k+1}T_k u=0
\quad\text{for all }u\in\dom(T_k).
\]
At the endpoints of the grading, the missing operator is understood to be
the zero operator.  We use the convention that the Hermitian inner product
is linear in its first argument and conjugate-linear in its second, so
that the adjoint is characterized by
$(Tu,v)=(u,T^{\ast}v)$.
Set
\[
\dom(Q):=\dom(T_k)\cap\dom(T_{k-1}^{\ast}).
\]
For \(u,v\in\dom(Q)\), define the Friedrichs form
\[
Q(u,v)
:=
(T_k u,T_k v)_{\mathrm{H}^{k+1}}+(T_{k-1}^{\ast}u,T_{k-1}^{\ast}v)_{\mathrm{H}^{k-1}}.
\]
Let \(\Box\) be the self-adjoint operator associated with the closed
nonnegative form \(Q\) by the Friedrichs representation theorem.  Then
\[
\dom(\Box)
=
\bigl\{u\in\dom(T_k)\cap\dom(T_{k-1}^{\ast}):
T_k u\in\dom(T_k^{\ast}),\quad
T_{k-1}^{\ast}u\in\dom(T_{k-1})\bigr\},
\]
and for every \(u\) in this domain,
\[
\Box u
=
T_k^{\ast}T_k u+T_{k-1}T_{k-1}^{\ast}u .
\]
\end{lemma}

\begin{proof}
\emph{ The form is closed and densely defined.}
The form \(Q\) is nonnegative.  It is closed with respect to the graph
norm
\[
\|u\|_Q^2:=\|u\|_{\mathrm{H}^k}^2+\|T_k u\|_{\mathrm{H}^{k+1}}^2
+\|T_{k-1}^{\ast}u\|_{\mathrm{H}^{k-1}}^2,
\]
because \(T_k\) and \(T_{k-1}^{\ast}\) are closed operators.

\emph{The adjoint-complex identity.}
Since \(T_kT_{k-1}=0\) in the strong sense, for every
\(v\in\dom(T_{k-1})\) we have \(T_{k-1}v\in\dom(T_k)\) and
\(T_kT_{k-1}v=0\). Hence, for \(u\in\dom(T_k^{\ast})\) and
\(v\in\dom(T_{k-1})\),
\[
(T_k^{\ast}u,T_{k-1}v)_{\mathrm{H}^{k}}
=
(u,T_kT_{k-1}v)_{\mathrm{H}^{k+1}}
=
0 .
\]
By the definition of the Hilbert adjoint, this shows that
\[
T_k^{\ast}u\in\dom(T_{k-1}^{\ast}),
\qquad
T_{k-1}^{\ast}T_k^{\ast}u=0
\qquad (u\in\dom(T_k^{\ast})),
\]
i.e.\ \(T_{k-1}^{\ast}T_k^{\ast}=0\) on \(\dom(T_k^{\ast})\).
Equivalently,
\[
\overline{\operatorname{Ran}T_k^{\ast}}\subset\ker T_{k-1}^{\ast}.
\]
We next show that \(\dom(Q)\) is dense in \(\mathrm{H}^k\).  For a closed, densely
defined operator \(T\), the domain \(\dom(T^{\ast}T)\) is dense in
\(\dom(T)\) with respect to the graph norm, and consequently
\(T(\dom(T^{\ast}T))\) is dense in \(\overline{\operatorname{Ran}T}\).
Applying this to \(T_{k-1}\) and to \(T_k^{\ast}\), and using the
orthogonal decomposition
\[
\mathrm{H}^k
=
\overline{\operatorname{Ran}T_{k-1}}
\oplus
\bigl(\ker T_k\cap\ker T_{k-1}^{\ast}\bigr)
\oplus
\overline{\operatorname{Ran}T_k^{\ast}},
\]
we obtain that \(\dom(Q)\) is dense in \(\mathrm{H}^k\).  (The decomposition
holds because \(\mathrm{H}^k=\overline{\operatorname{Ran}T_{k-1}}\oplus\ker T_{k-1}^{\ast}\),
\(\mathrm{H}^k=\ker T_k\oplus\overline{\operatorname{Ran}T_k^{\ast}}\), and
\(\overline{\operatorname{Ran}T_k^{\ast}}\subset\ker T_{k-1}^{\ast}\) by the
adjoint-complex identity, so that
\(\ker T_{k-1}^{\ast}=(\ker T_k\cap\ker T_{k-1}^{\ast})\oplus\overline{\operatorname{Ran}T_k^{\ast}}\).)
Indeed, the first
summand is approximated by
\(T_{k-1}(\dom(T_{k-1}^{\ast}T_{k-1}))\subset\dom(Q)\), the middle
summand is contained in \(\dom(Q)\), and the third summand is approximated
by \(T_k^{\ast}(\dom(T_kT_k^{\ast}))\subset\dom(Q)\).  Thus the
Friedrichs representation theorem applies and defines a unique
self-adjoint operator \(\Box\) with
\[
Q(u,v)=(\Box u,v)_\mathrm{H}
\qquad
(u\in\dom(\Box),\ v\in\dom(Q)).
\]

\emph{ The forward domain inclusion $\subseteq$ direction.}
Let \(u\in\dom(\Box)\) and put \(w:=\Box u\in \mathrm{H}^k\).  Then
\(u\in\dom(Q)\) and

\[
Q(u,v)=(w,v)_{\mathrm{H}^k}
\qquad\text{for all }v\in\dom(Q).
\]
We first prove \(T_k u\in\dom(T_k^{\ast})\).  Let \(\mathrm{P}_{\!A}\) be the
orthogonal projection of \(\mathrm{H}^k\) onto \(\ker T_{k-1}^{\ast}\).  Since
\[
(\ker T_{k-1}^{\ast})^\perp
=
\overline{\operatorname{Ran}T_{k-1}},
\]
for every \(v\in\dom(T_k)\) we have
\[
(I-\mathrm{P}_{\!A})v\in\overline{\operatorname{Ran}T_{k-1}}
\subset\ker T_k\subset\dom(T_k).
\]
The inclusion \(\overline{\operatorname{Ran}T_{k-1}}\subset\ker T_k\)
follows from \(T_kT_{k-1}=0\), the domain inclusion
\(\operatorname{Ran}T_{k-1}\subset\dom(T_k)\), and the closedness of
\(T_k\).  Hence
\[
\mathrm{P}_{\!A}v=v-(I-\mathrm{P}_{\!A})v\in\dom(T_k),
\qquad
T_k\mathrm{P}_{\!A}v=T_kv.
\]
Moreover \(\mathrm{P}_{\!A}v\in\ker T_{k-1}^{\ast}\subset\dom(T_{k-1}^{\ast})\), so
\(\mathrm{P}_{\!A}v\in\dom(Q)\).  Using \(T_{k-1}^{\ast}\mathrm{P}_{\!A}v=0\), we obtain
\[
(T_k u,T_k v)_\mathrm{H}
=
(T_k u,T_k \mathrm{P}_{\!A}v)_\mathrm{H}
=
Q(u,\mathrm{P}_{\!A}v)
=
(w,\mathrm{P}_{\!A}v)_\mathrm{H}
=
(\mathrm{P}_{\!A}w,v)_\mathrm{H} .
\]
Thus the conjugate-linear functional
\(v\mapsto(T_k u,T_k v)_{\mathrm{H}^{k+1}}\) is bounded on \(\dom(T_k)\).
By the definition of the Hilbert adjoint,
\[
T_k u\in\dom(T_k^{\ast}),
\qquad
T_k^{\ast}T_k u=\mathrm{P}_{\!A}w .
\]

Next we prove \(T_{k-1}^{\ast}u\in\dom(T_{k-1})\).  Let \(\mathrm{P}_T\) be the
orthogonal projection of \(\mathrm{H}^k\) onto \(\ker T_k\).  Since
\[
(\ker T_k)^\perp=\overline{\operatorname{Ran}T_k^{\ast}},
\]
for every \(v\in\dom(T_{k-1}^{\ast})\) we have
\[
(I-\mathrm{P}_T)v\in\overline{\operatorname{Ran}T_k^{\ast}}
\subset\ker T_{k-1}^{\ast}
\subset\dom(T_{k-1}^{\ast}).
\]
By the adjoint-complex identity proved above,
\(\overline{\operatorname{Ran}T_k^{\ast}}\subset\ker T_{k-1}^{\ast}\). Hence
\[
\mathrm{P}_Tv=v-(I-\mathrm{P}_T)v\in\dom(T_{k-1}^{\ast}),
\qquad
T_{k-1}^{\ast}\mathrm{P}_Tv=T_{k-1}^{\ast}v.
\]
Moreover \(\mathrm{P}_Tv\in\ker T_k\subset\dom(T_k)\), and \(T_k\mathrm{P}_Tv=0\).  Thus
\(\mathrm{P}_Tv\in\dom(Q)\).  Testing the form identity with \(\mathrm{P}_Tv\) gives
\[
(T_{k-1}^{\ast}u,T_{k-1}^{\ast}v)_\mathrm{H}
=
Q(u,\mathrm{P}_Tv)
=
(w,\mathrm{P}_Tv)_\mathrm{H}
=
(\mathrm{P}_Tw,v)_\mathrm{H} .
\]
Since \(T_{k-1}^{\ast\ast}=T_{k-1}\), the definition of the adjoint yields
\[
T_{k-1}^{\ast}u\in\dom(T_{k-1}),
\qquad
T_{k-1}T_{k-1}^{\ast}u=\mathrm{P}_Tw .
\]
Therefore \(u\) belongs to the set displayed in the lemma.  Moreover,
since \(\dom(Q)\) is dense and
\[
Q(u,v)
=
(\mathrm{P}_{\!A}w+\mathrm{P}_Tw,v)_\mathrm{H}
\qquad
(v\in\dom(Q)),
\]
we have \(w=\mathrm{P}_{\!A}w+\mathrm{P}_Tw\).  Consequently
\[
\Box u
=
T_k^{\ast}T_k u+T_{k-1}T_{k-1}^{\ast}u .
\]

\emph{ The reverse domain inclusion $\supseteq$ direction.}
Conversely, suppose that
\[
u\in\dom(T_k)\cap\dom(T_{k-1}^{\ast}),
\qquad
T_k u\in\dom(T_k^{\ast}),
\qquad
T_{k-1}^{\ast}u\in\dom(T_{k-1}).
\]
Then for every \(v\in\dom(Q)\),
\[
\begin{aligned}
Q(u,v)
&=
(T_k u,T_k v)_\mathrm{H}+(T_{k-1}^{\ast}u,T_{k-1}^{\ast}v)_\mathrm{H}\\
&=
(T_k^{\ast}T_k u,v)_\mathrm{H}+(T_{k-1}T_{k-1}^{\ast}u,v)_\mathrm{H}\\
&=
\bigl(T_k^{\ast}T_k u+T_{k-1}T_{k-1}^{\ast}u,v\bigr)_\mathrm{H} .
\end{aligned}
\]
Thus \(u\in\dom(\Box)\), and
\[
\Box u=T_k^{\ast}T_k u+T_{k-1}T_{k-1}^{\ast}u .
\]
\end{proof}

\begin{definition}[Twisted \(\overline{\partial}_f\)-Neumann problem]
\label{def:twisted-neumann}
Let \(0\le k\le 2n\), \(\alpha\in L^2_{\Lambda^k}(\Omega)\), and put
\[
\dom(Q_f):=\dom(\overline{\partial}_f)\cap\dom(\overline{\partial}_f^{\ast})
\subset L^2_{\Lambda^k}(\Omega),
\]
\[
Q_f(u,v):=(\overline{\partial}_f^{\ast} u,\overline{\partial}_f^{\ast} v)+(\overline{\partial}_fu,\overline{\partial}_fv),
\qquad u,v\in\dom(Q_f).
\]
The \emph{weak twisted \(\overline{\partial}_f\)-Neumann problem} is to find
\(u\in\dom(Q_f)\) such that
\[
Q_f(u,v)=(\alpha,v)
\qquad\text{for all }v\in\dom(Q_f).
\]
Since \(\overline{\partial}_f\) and \(\overline{\partial}_f^{\ast}\) are closed, \(Q_f\) is a closed
nonnegative form.  It is also densely defined: every compactly supported
smooth \(k\)-form \(u\in C_{c,\Lambda^k}^\infty(\Omega)\) has
\(\operatorname{supp}u\Subset\Omega\), hence belongs to
\(\dom(\overline{\partial}_f)\cap\dom(\overline{\partial}_f^{\ast})\)
and therefore to \(\dom(Q_f)\), and such forms are dense in
\(L^2_{\Lambda^k}(\Omega)\).  Thus \(Q_f\) is a closed, densely defined,
nonnegative form.  By the Friedrichs representation theorem, there is a unique
self-adjoint operator associated with \(Q_f\), denoted \(\Box_f\), defined by
\[
\dom(\Box_f):=
\bigl\{u\in\dom(Q_f):\ \exists\,w\in L^2_{\Lambda^k}(\Omega)
\text{ with }
Q_f(u,v)=(w,v)\ \ \forall v\in\dom(Q_f)
\bigr\},
\]
with \(\Box_fu:=w\). 
Both \(\overline{\partial}_f\) and \(\overline{\partial}_f^{\ast}\) are closed and densely defined: \(\overline{\partial}_{f,k}\) is the maximal realization of a differential operator with smooth coefficients, hence closed (if \(u_j\to u\) and \(\overline{\partial}_fu_j\to g\) in \(L^2\), then \(\overline{\partial}_fu=g\) in the sense of distributions), and it contains \(C^\infty_{c,\Lambda^k}(\Omega)\); \(\overline{\partial}_{f,k}^{\ast}\) is a Hilbert adjoint, hence closed, and contains \(C^\infty_{c,\Lambda^k}(\Omega)\) by Lemma~\ref{lem:domains}(a),(c).  Moreover \(\overline{\partial}_{f,k+1}\overline{\partial}_{f,k}=0\) in the strong sense of Lemma~\ref{lem:abstract-Hilbert-complex}: for \(u\in\dom(\overline{\partial}_{f,k})\), the distribution \(\overline{\partial}_f(\overline{\partial}_fu)\) vanishes, hence lies in \(L^2\), so \(\overline{\partial}_fu\in\dom(\overline{\partial}_{f,k+1})\) by maximality of the realization and \(\overline{\partial}_{f,k+1}\overline{\partial}_{f,k}u=0\).  Thus the twisted operators form a Hilbert complex. By Lemma~\ref{lem:abstract-Hilbert-complex}, we obtain
\[
\dom(\Box_f)
=
\bigl\{u\in\dom(\overline{\partial}_f)\cap\dom(\overline{\partial}_f^{\ast}):
\overline{\partial}_f^{\ast} u\in\dom(\overline{\partial}_f),\quad
\overline{\partial}_fu\in\dom(\overline{\partial}_f^{\ast})
\bigr\},
\]
and on this domain
\[
\Box_fu
=
\overline{\partial}_f\overline{\partial}_f^{\ast} u+\overline{\partial}_f^{\ast} \overline{\partial}_fu
\qquad\text{in }L^2_{\Lambda^k}(\Omega).
\]
We refer to \(\Box_f\) as the \emph{twisted complex Laplacian}.
The \emph{strong form} of the problem is: find
\(u\in\dom(\Box_f)\) such that
\[
\Box_fu=\alpha
\qquad\text{in }L^2_{\Lambda^k}(\Omega).
\]
It is equivalent to the weak formulation, because
\(Q_f(u,v)=(\alpha,v)\) for every \(v\in\dom(Q_f)\) says precisely that
the functional \(v\mapsto Q_f(u,v)\) is represented by
\(\alpha\in L^2_{\Lambda^k}(\Omega)\), which is exactly the defining
condition for \(u\in\dom(\Box_f)\) and \(\Box_fu=\alpha\).
The weak Neumann boundary condition is already encoded in
\(u\in\dom(\overline{\partial}_f^{\ast})\), hence in \(\dom(Q_f)\).
\end{definition}

\begin{remark}[Conventions on \(H^\epsilon\)-norms]
\label{rem:Heps-norms}
Throughout, for an open set \(U\subset\mathbb C^n\), the symbol
\(\|F\|_{H^\epsilon(U)}\) denotes the intrinsic Gagliardo norm
\[
\|F\|_{H^\epsilon(U)}^2
=
\|F\|_{L^2(U)}^2
+
\iint_{U\times U}
\frac{|F(x)-F(y)|^2}{|x-y|^{2n+2\epsilon}}\,dx\,dy .
\]
For \(U=\mathbb C^n=\mathbb R^{2n}\) this norm is equivalent, with
constants depending only on \(n\) and \(\epsilon\), to the Fourier norm
\(\bigl(\int(1+|\xi|^2)^{\epsilon}|\widehat F|^2\,d\xi\bigr)^{1/2}\)
(cf.\ \cite[7.62 and Theorem~7.47]{AdamFou}): by Plancherel,
\(\iint|F(x+h)-F(x)|^2|h|^{-2n-2\epsilon}\,dx\,dh
=\int|\widehat F(\xi)|^2\int|e^{ih\cdot\xi}-1|^2|h|^{-2n-2\epsilon}\,dh\,d\xi
=c_{2n,\epsilon}\int|\xi|^{2\epsilon}|\widehat F(\xi)|^2\,d\xi\) with
\(0<c_{2n,\epsilon}<\infty\) for \(0<\epsilon<1\).
For any open set \(U\subset\mathbb C^n\), we also use the quotient
(restriction) norm
\[
\|F\|_{H^\epsilon(U)}^{\mathrm{quot}}
:=
\inf\bigl\{
\|G\|_{H^\epsilon(\mathbb C^n)}:\ G|_U=F
\bigr\}.
\]
The direction needed below is the elementary inequality
$\|F\|_{H^\epsilon(U)}
\le
C(n,\epsilon)\,\|F\|_{H^\epsilon(U)}^{\mathrm{quot}}$,
which holds for every open set \(U\). If the \(H^\epsilon(\mathbb C^n)\)-norm
used in the definition of the quotient norm is the Gagliardo norm, then the
inequality follows by restricting the Gagliardo double integral over
\(\mathbb C^n\times\mathbb C^n\) to \(U\times U\) and using
\(\|F\|_{L^2(U)}\le\|G\|_{L^2(\mathbb C^n)}\) for any extension \(G\);
in this case one may take \(C(n,\epsilon)=1\). If instead another norm
equivalent to the Gagliardo norm is used on \(\mathbb C^n\), the same
argument gives the inequality with \(C(n,\epsilon)\) equal to the equivalence
constant between the two norms, which depends only on \(n\) and \(\epsilon\).

Two further elementary facts are used.  First, for
\(0<\epsilon\le\epsilon^{\prime}<1\) and every open \(U\subset\mathbb C^n\),
\begin{equation}\label{eq:Heps-monotone}
\|F\|_{H^\epsilon(U)}^2\le\bigl(1+4A_{n,\epsilon}\bigr)\|F\|_{H^{\epsilon^{\prime}}(U)}^2,
\qquad
A_{n,\epsilon}:=\int_{|h|\ge1}|h|^{-2n-2\epsilon}\,dh<\infty ;
\end{equation}
indeed, on \(\{|x-y|<1\}\) one has
\(|x-y|^{-2n-2\epsilon}\le|x-y|^{-2n-2\epsilon^{\prime}}\), while on
\(\{|x-y|\ge1\}\) one uses \(|F(x)-F(y)|^2\le2|F(x)|^2+2|F(y)|^2\).
Second, if \(F\in H^1(\mathbb C^n)\) (for instance the zero extension of
a compactly supported \(H^1\)-form on an open subset of \(\mathbb C^n\)),
then for every open \(U\subset\mathbb C^n\) and \(0<\epsilon<1\),
\begin{equation}\label{eq:H1-Heps}
\|F\|_{H^\epsilon(U)}^2\le\|F\|_{H^\epsilon(\mathbb C^n)}^2
\le\bigl(1+c_{2n,\epsilon}\bigr)\bigl(\|F\|_{L^2(\mathbb C^n)}^2+\|\nabla F\|_{L^2(\mathbb C^n)}^2\bigr),
\end{equation}
by restriction of the double integral and the Plancherel identity above
together with \(|\xi|^{2\epsilon}\le1+|\xi|^2\).
\end{remark}

For \(0<\epsilon<1\) and open \(U\subset\mathbb C^n=\mathbb R^{2n}\), write
\[
[g]_{H^\epsilon(U)}^2
:=
\iint_{U\times U}
\frac{|g(x)-g(y)|^2}{|x-y|^{2n+2\epsilon}}\,dx\,dy
\]
for the Gagliardo seminorm.

\emph{Straube estimate neighborhoods.}
For each \(P\in\partial\Omega\), the boundary geometry hypothesis gives
a neighborhood \(U_P\) and a bounded plurisubharmonic function
\(\lambda_P\) satisfying \eqref{eq:boundary-hyp} with exponent
\(\epsilon_P\in(0,1/2)\).  By \cite[Theorem~1]{Straube1997} (whose hypotheses
are: \(\Omega\) bounded pseudoconvex, \(\partial\Omega\) a Lipschitz graph near
\(P\), and \eqref{eq:boundary-hyp} with \(0\le\epsilon_P<1/2\)), there exist a
neighborhood \(V_P\) of \(P\) and a constant \(C_P>0\) such that the estimate below
holds on \(V_P\cap\Omega\); it persists on any smaller neighborhood, since
the restriction norm decreases with the set.  We choose
\(U_P^{\mathrm{est}}\Subset U_P\cap V_P\) with \(P\in U_P^{\mathrm{est}}\).
Then for every \(1\le q\le n\) and every
\((0,q)\)-form
\[
\varphi\in L^2_{(0,q)}(\Omega)\cap\dom(\overline{\partial})
\cap\dom(\overline{\partial}^{\ast}),
\]
the following boundary-local subelliptic estimate holds:
\begin{equation}\label{eq:straube-local-est}
\|\varphi\|_{H^{\epsilon_P}(U_P^{\mathrm{est}}\cap\Omega)}
\le C_P\Bigl(
\|\overline{\partial}\varphi\|_{L^2(\Omega)}
+\|\overline{\partial}^{\ast}\varphi\|_{L^2(\Omega)}
\Bigr).
\end{equation}
Here the left-hand norm is the intrinsic Gagliardo norm; the
quotient-norm form of \cite[Theorem~1]{Straube1997} is converted to the
intrinsic norm by Remark~\ref{rem:Heps-norms}, at the cost of
enlarging \(C_P\).  Since \(\partial\Omega\) is compact, there exist
finitely many \(P_1,\dots,P_M\in\partial\Omega\) such that
\[
\partial\Omega\subset\bigcup_{\ell=1}^M U_\ell^{\mathrm{est}},
\qquad
U_\ell^{\mathrm{est}}:=U_{P_\ell}^{\mathrm{est}}\Subset U_{P_\ell}.
\]
Set
\[
\epsilon:=\min_{1\le\ell\le M}\epsilon_{P_\ell}>0 .
\]
After enlarging the constants \(C_\ell:=C_{P_\ell}\), estimate
\eqref{eq:straube-local-est} holds with the uniform exponent
\(\epsilon\) on each \(U_\ell^{\mathrm{est}}\cap\Omega\), because
\(0<\epsilon\le\epsilon_{P_\ell}\) and, by
\eqref{eq:Heps-monotone} in Remark~\ref{rem:Heps-norms},
\(\|\varphi\|_{H^{\epsilon}(U\cap\Omega)}\le C(n,\epsilon)\|\varphi\|_{H^{\epsilon_{P_\ell}}(U\cap\Omega)}\)
for every open \(U\).  Moreover, since
\(U_\ell^{\mathrm{est}}\subset U_{P_\ell}\), the psh functions
\(\lambda_{P_\ell}\) satisfy
\[
dd^c\lambda_{P_\ell}
\ge \operatorname{const.}\delta^{-2\epsilon}
dd^c|z|^2
\qquad\text{on }U_\ell^{\mathrm{est}}\cap\Omega .
\]
All local subelliptic estimates below are understood with this uniform
exponent \(\epsilon\) and with these estimate neighborhoods
\(U_\ell^{\mathrm{est}}\); the constants \(C_\ell\) are taken as the
maximum over \(1\le q\le n\).

\begin{remark}
The inequality in \eqref{eq:boundary-hyp} is an inequality between
\((1,1)\)-currents on the open set \(U_P\cap\Omega\). The boundary
\(\partial\Omega\) enters only through the Lipschitz graph assumption and
through the distance function \(\delta\); it is not itself treated as a
current.
\end{remark}

\emph{Componentwise reduction in the holomorphic degree.}
We use the elementary reduction of \((p,q)\)-forms to \((0,q)\)-forms
available in flat \(\mathbb C^n\) (for the general \((p,q)\) setting of
the \(\overline\partial\)-Neumann problem see \cite{FollandKohn1972}).
On \(\mathbb C^n\) with the flat metric, the bundle \(\Lambda^{p,0}\) is
holomorphically and metrically trivial, with the constant orthonormal
coframe \(\{dz_I:|I|=p\}\).
Writing a \((p,q)\)-form as
\[
u=\sum_{|I|=p} dz_I\wedge u_I,
\qquad u_I\in\Lambda^{0,q},
\]
the operators in play act only on the antiholomorphic factor:
\[
\overline\partial u=(-1)^p\sum_{|I|=p}dz_I\wedge\overline\partial u_I,
\quad
\overline\partial^{\ast} u=(-1)^p\sum_{|I|=p}dz_I\wedge\overline\partial^{\ast} u_I,
\]
and, wherever \(\overline{\partial}\rho\) is defined,
$\iota_{(\overline{\partial}\rho)^\sharp}u=(-1)^p\sum_{|I|=p}dz_I\wedge\iota_{(\overline{\partial}\rho)^\sharp}u_I$,
the last identity because \((\overline{\partial}\rho)^\sharp\) is a \((0,1)\)-vector and
therefore annihilates the holomorphic factor \(dz_I\). The precise
statement used below, including the splitting of the Hilbert-space
domains, is Lemma~\ref{lem:componentwise}.

\begin{lemma}[Componentwise reduction in the holomorphic degree]
\label{lem:componentwise}
Let \(\Omega\subset\mathbb C^n\) be a bounded domain, \(0\le p\le n\),
\(0\le q\le n\), and let \(\mathcal I_p\) be the set of strictly
increasing \(p\)-tuples.  Define
\[
\Phi_p:\bigoplus_{I\in\mathcal I_p}L^2_{(0,q)}(\Omega)\to L^2_{(p,q)}(\Omega),
\qquad
(v_I)_I\mapsto\sum_{I}dz_I\wedge v_I .
\]
\begin{enumerate}[label=\textup{(\alph*)}]
\item \(\Phi_p\) is unitary, and pointwise
\(|\Phi_p(v_I)_I|^2=\sum_I|v_I|^2\).
\item \(u=\Phi_p(u_I)_I\in\dom(\overline\partial)\) in \(L^2_{(p,q)}(\Omega)\) if and
only if \(u_I\in\dom(\overline\partial)\) in \(L^2_{(0,q)}(\Omega)\) for every \(I\);
then \(\overline\partial u=(-1)^p\Phi_p(\overline\partial u_I)_I\).
\item For \(q\ge1\), \(u=\Phi_p(u_I)_I\in\dom(\overline\partial^{\ast})\) in
\(L^2_{(p,q)}(\Omega)\) if and only if \(u_I\in\dom(\overline\partial^{\ast})\)
in \(L^2_{(0,q)}(\Omega)\) for every \(I\); then
\(\overline\partial^{\ast}u=(-1)^p\Phi_p(\overline\partial^{\ast}u_I)_I\).
\item For every open \(U\subset\mathbb C^n\) and \(0<\epsilon<1\),
\(\|u\|_{H^\epsilon(U\cap\Omega)}^2=\sum_I\|u_I\|_{H^\epsilon(U\cap\Omega)}^2\)
for the intrinsic Gagliardo norm.
\end{enumerate}
Consequently, let \(1\le q\le n\) and \(U\subset\mathbb C^n\) be open, and assume that an estimate of the form \eqref{eq:straube-local-est},
\(\|\varphi\|_{H^\epsilon(U\cap\Omega)}\le C(\|\overline\partial\varphi\|_{L^2(\Omega)}+\|\overline\partial^{\ast}\varphi\|_{L^2(\Omega)})\), holds
for all \(\varphi\in L^2_{(0,q)}(\Omega)\cap\dom(\overline\partial)\cap\dom(\overline\partial^{\ast})\).
Then
\(\|v\|_{H^\epsilon(U\cap\Omega)}\le\sqrt2\,C(\|\overline\partial v\|_{L^2(\Omega)}+\|\overline\partial^{\ast}v\|_{L^2(\Omega)})\)
for all \(v\in L^2_{(p,q)}(\Omega)\cap\dom(\overline\partial)\cap\dom(\overline\partial^{\ast})\),
with \(\sqrt2\,C\) independent of \(p\).
\end{lemma}

\begin{proof}
(a) Write \(v_I=\sum_J v_{IJ}\,d\overline z_J\).  Then
\(\Phi_p(v_I)_I=\sum_{I,J}v_{IJ}\,dz_I\wedge d\overline z_J\), and since the
coframe \(\{dz_I\wedge d\overline z_J\}\) is orthonormal,
\(|\Phi_p(v_I)_I|^2=\sum_{I,J}|v_{IJ}|^2=\sum_I|v_I|^2\) pointwise.
Every \((p,q)\)-form has a unique expansion of this type, so \(\Phi_p\)
is a surjective isometry, hence unitary.

(b) The coefficients of \(dz_I\) are constant, so
\(\overline\partial(dz_I\wedge v_I)=(-1)^p\,dz_I\wedge\overline\partial v_I\) as an identity of
constant-coefficient differential expressions in the \(v_{IJ}\); it
therefore holds in the sense of distributions.  Thus
\(\overline\partial u=(-1)^p\Phi_p(\overline\partial u_I)_I\) as distributions, and by (a) in
bidegree \((p,q+1)\), \(\overline\partial u\in L^2\) if and only if every
\(\overline\partial u_I\in L^2\).  Since \(\dom(\overline\partial)\) is the maximal domain, (b)
follows.

(c) By (b) in bidegree \((p,q-1)\),
\(\overline\partial_{(p,q-1)}=\Phi_p\,\mathcal{T}\,\Phi_p^{-1}\) with
\(\mathcal{T}:=(-1)^p\bigoplus_I\overline\partial_{(0,q-1)}\) and
\(\dom(\mathcal{T})=\bigoplus_I\dom(\overline\partial)\).  As \(\Phi_p\) is unitary,
\(\overline\partial^{\ast}=\Phi_p\,\mathcal{T}^{\ast}\,\Phi_p^{-1}\).  For a finite
orthogonal direct sum of densely defined operators \(A_I\), one has
\((\bigoplus_IA_I)^{\ast}=\bigoplus_IA_I^{\ast}\) with domain
\(\bigoplus_I\dom(A_I^{\ast})\): the functional
\((v_I)_I\mapsto\sum_I(w_I,A_Iv_I)\) is bounded on
\(\bigoplus_I\dom(A_I)\) if and only if each
\(v_I\mapsto(w_I,A_Iv_I)\) is bounded on \(\dom(A_I)\), as one sees by
testing with families having a single nonzero entry.  This gives (c).

(d) Applying (a) to \(u(x)\) and to \(u(x)-u(y)\) gives
\(|u(x)|^2=\sum_I|u_I(x)|^2\) and
\(|u(x)-u(y)|^2=\sum_I|u_I(x)-u_I(y)|^2\); integrating over
\(U\cap\Omega\) and \((U\cap\Omega)^2\) against the Gagliardo kernel yields (d).

For the consequence, let \(v=\Phi_p(v_I)_I\).  By (b) and (c), each
\(v_I\in L^2_{(0,q)}(\Omega)\cap\dom(\overline\partial)\cap\dom(\overline\partial^{\ast})\).
Put \(a_I=\|\overline\partial v_I\|\), \(b_I=\|\overline\partial^{\ast}v_I\|\).  By
(d), the hypothesis, \((a_I+b_I)^2\le2(a_I^2+b_I^2)\), and (a) together
with (b), (c),
\[
\|v\|_{H^\epsilon(U\cap\Omega)}^2
\le C^2\sum_I(a_I+b_I)^2
\le2C^2\bigl(\|\overline\partial v\|^2+\|\overline\partial^{\ast}v\|^2\bigr)
\le2C^2\bigl(\|\overline\partial v\|+\|\overline\partial^{\ast}v\|\bigr)^2 .
\]
When the action is evident from the context, we suppress the subscripts on the operators.
\end{proof}

\emph{Action on \(\Lambda^{p,q}\) and the bidegree-free bridge.}
Let \(u\in\Lambda^{p,q}\) with \(0\le p\le n\) and  \(1\le q\le n\). The four summands of \(\overline{\partial}_fu\) and \(\overline{\partial}_f^{\ast}u\),
\[
\overline\partial u\in\Lambda^{p,q+1},\quad
\partial f\wedge u\in\Lambda^{p+1,q},\quad
\overline\partial^{\ast} u\in\Lambda^{p,q-1},\quad
\iota_{(\partial f)^\sharp}u\in\Lambda^{p-1,q},
\]
lie in four \emph{distinct} bidegrees. (For \(p=0\) the last summand
\(\iota_{(\partial f)^\sharp}u\) vanishes identically, since
\((\partial f)^\sharp\) is of type \((1,0)\) and \(u\) has no holomorphic
degree; this recovers the Straube case \(\overline{\partial}_f^{\ast}u=\overline\partial^{\ast} u\).) By pointwise
orthogonality of distinct bidegrees,
\begin{equation*}
\|\overline{\partial}_fu\|^2=\|\overline\partial u\|^2+\|\partial f\wedge u\|^2,
\qquad
\|\overline{\partial}_f^{\ast}u\|^2=\|\overline\partial^{\ast} u\|^2+\|\iota_{(\partial f)^\sharp}u\|^2 .
\end{equation*}
Both identities are exact; the cross terms vanish because the paired forms
have different bidegrees. Dropping the two nonnegative zero-order
contributions gives the \emph{bidegree-free bridge inequality}
\begin{equation}
\label{eq:bridge}
\|\overline\partial u\|^2+\|\overline\partial^{\ast} u\|^2
\le
\|\overline{\partial}_fu\|^2+\|\overline{\partial}_f^{\ast}u\|^2 ,
\qquad u\in L^2_{(p,q)}(\Omega)\cap\dom(\overline{\partial}_f)\cap\dom(\overline{\partial}_f^{\ast}).
\end{equation}
For such \(u\), Lemma~\ref{lem:domains}(a) gives \(\overline\partial u,\overline\partial^{\ast}u\in L^2\),
and the pointwise orthogonality of the four distinct bidegrees yields
\(|\overline\partial u|^2+|\overline\partial^{\ast}u|^2\le|\overline{\partial}_fu|^2+|\overline{\partial}_f^{\ast}u|^2\)
almost everywhere; integrating gives \eqref{eq:bridge}, and integrating
against any nonnegative weight, e.g.\ \(\theta^2\), gives its weighted form
\begin{equation*}
\int_\Omega\theta^2\bigl(|\overline\partial u|^2+|\overline\partial^{\ast} u|^2\bigr)\,dV
\le
\int_\Omega\theta^2\bigl(|\overline{\partial}_fu|^2+|\overline{\partial}_f^{\ast}u|^2\bigr)\,dV .
\end{equation*}
This is the bridge that transports the untwisted \(\overline\partial\)-Neumann estimates
to the twisted operators.  It holds for every bidegree \((p,q)\): the
twisted zero-order terms appear on the dominating side and are
discarded, so no term \(C\|u\|\) is produced and no property of
\(\partial f\) beyond boundedness is used.  It is used only for forms of a
fixed bidegree \((p,q)\).

\section{The \texorpdfstring{\(L^2\)-Version}{L2 Version} of the Multiplier and Collar-Piece Estimates}

\begin{lemma}[Multiplication by real \(W^{1,\infty}\)-multipliers on a fixed bidegree]
\label{lem:multiplier-closure-L2}
Let \(k=p+q\). 

Let \(\overline{\partial}_f\) be the maximal \(L^2\)-realization
\[
\overline{\partial}_f:\dom(\overline{\partial}_f)\subset L^2_{\Lambda^k}(\Omega)\to L^2_{\Lambda^{k+1}}(\Omega),
\]
and let \(\overline{\partial}_f^{\ast}\) be its Hilbert adjoint in the graded-adjoint convention, so that
\[
\overline{\partial}_f^{\ast}:\dom(\overline{\partial}_f^{\ast})\subset L^2_{\Lambda^k}(\Omega)\to L^2_{\Lambda^{k-1}}(\Omega)
\]
is the Hilbert adjoint of \(\overline{\partial}_f:L^2_{\Lambda^{k-1}}\to L^2_{\Lambda^k}\).
For a \((p,q)\)-form \(u\in L^2_{(p,q)}(\Omega)\), the conditions
\(u\in\dom(\overline{\partial}_f)\) and \(u\in\dom(\overline{\partial}_f^{\ast})\) are understood in this total-degree
sense; equivalently,
\[
\overline{\partial}_fu=\overline\partial u+\partial f\wedge u\in L^2_{(p,q+1)}(\Omega)\oplus L^2_{(p+1,q)}(\Omega),
\qquad
\overline{\partial}_f^{\ast}u=\overline\partial^{\ast}u+\iota_{(\partial f)^\sharp}u
\in L^2_{(p,q-1)}(\Omega)\oplus L^2_{(p-1,q)}(\Omega).
\]
Let
\[
u\in L^2_{(p,q)}(\Omega)\cap\dom(\overline{\partial}_f)\cap\dom(\overline{\partial}_f^{\ast}),
\qquad
\sigma\in W^{1,\infty}(\Omega;\mathbb R).
\]
We use the convention
\[
\Lambda^{r,s}=0
\qquad\text{whenever }r<0,\ r>n,\ s<0,\ \text{or }s>n.
\]
Then
\[
\sigma u\in L^2_{(p,q)}(\Omega)\cap\dom(\overline{\partial}_f)\cap\dom(\overline{\partial}_f^{\ast}),
\]
and hence \(\sigma u\) also satisfies the weak \(\overline{\partial}_f\)-Neumann boundary
condition. Moreover, the weak product rules
\begin{equation}\label{eq:weak-product}
\overline{\partial}_f(\sigma u)=\sigma \overline{\partial}_fu+\overline{\partial}\sigma\wedge u,
\qquad
\overline{\partial}_f^{\ast}(\sigma u)=\sigma \overline{\partial}_f^{\ast}u-\iota_{(\overline{\partial}\sigma)^\sharp}u
\end{equation}
hold in \(L^2\).
\end{lemma}

\begin{proof}
Since \(u\in L^2_{(p,q)}(\Omega)\) and \(\sigma\in W^{1,\infty}(\Omega;\mathbb R)\subset L^\infty(\Omega)\), we have
\[
\sigma u\in L^2_{(p,q)}(\Omega).
\]

We first show that \(u\in\dom(\overline{\partial})\cap\dom(\overline{\partial}^{\ast})\).
On the fixed bidegree \(\Lambda^{p,q}\),
\[
\overline{\partial}_f=\overline{\partial}+\partial f\wedge,
\qquad
\overline{\partial}_f^{\ast}=\overline{\partial}^{\ast}+\iota_{(\partial f)^\sharp},
\]
and the two summands in each operator take values in distinct bidegrees:
\[
\overline{\partial} u\in\Lambda^{p,q+1},
\quad
\partial f\wedge u\in\Lambda^{p+1,q},
\qquad
\overline{\partial}^{\ast}u\in\Lambda^{p,q-1},
\quad
\iota_{(\partial f)^\sharp}u\in\Lambda^{p-1,q}.
\]
Since \(u\in\dom(\overline{\partial}_f)\) and \(u\in\dom(\overline{\partial}_f^{\ast})\), we have
\(\overline{\partial}_fu\in L^2\) and \(\overline{\partial}_f^{\ast}u\in L^2\), while
\(\partial f\wedge u\) and \(\iota_{(\partial f)^\sharp}u\) are bounded
zero-order terms and hence belong to \(L^2\) as well. Therefore, the formal expressions satisfy
 \[
 \overline{\partial} u=\overline{\partial}_fu-\partial f\wedge u\in L^2,
 \qquad
 \overline{\partial}^{\ast}u=\overline{\partial}_f^{\ast}u-\iota_{(\partial f)^\sharp}u\in L^2 .
 \]
 Moreover, by the domain equality,
 \[
 \dom(\overline{\partial}_f)=\dom(\overline{\partial}),
 \qquad
 \dom(\overline{\partial}_f^{\ast})=\dom(\overline{\partial}^{\ast}),
 \]
 and since \(u\in\dom(\overline{\partial}_f)\cap\dom(\overline{\partial}_f^{\ast})\), we obtain
 \[
 u\in\dom(\overline{\partial})\cap\dom(\overline{\partial}^{\ast}),
 \]
 with \(\overline{\partial} u\) and \(\overline{\partial}^{\ast}u\) given by the expressions above.

Next we prove that \(\sigma u\in\dom(\overline{\partial})\). Because
\(\sigma\in W^{1,\infty}(\Omega)\), the Leibniz identity for
\(\overline{\partial}\) holds in the sense of distributions:
\[
\overline{\partial}(\sigma u)=\sigma\,\overline{\partial} u+\overline{\partial}\sigma\wedge u.
\]
Indeed, let \(\vartheta\) denote the formal adjoint of \(\overline{\partial}\)
and let \(\varphi\) be a smooth compactly supported test form in \(\Omega\).
Then \(\sigma\varphi\) is a compactly supported \(W^{1,\infty}\)-form, and
let \(J_\eta\) be a standard mollifier and set
\((\sigma\varphi)_\eta:=J_\eta \ast(\sigma\varphi)\in C^\infty_c(\Omega)\).
Then \((\sigma\varphi)_\eta\to\sigma\varphi\) in \(L^2\), and, since
\(\vartheta\) is a constant-coefficient differential expression and hence
commutes with convolution,
\[
\vartheta\bigl((\sigma\varphi)_\eta\bigr)
=J_\eta \ast \bigl(\vartheta(\sigma\varphi)\bigr)
\to\vartheta(\sigma\varphi)
=\sigma\vartheta\varphi-\iota_{(\overline{\partial}\sigma)^\sharp}\varphi
\qquad\text{in }L^2 .
\]
Passing to the limit in
\[
\bigl(u,\vartheta((\sigma\varphi)_\eta)\bigr)
=\bigl(\overline{\partial} u,(\sigma\varphi)_\eta\bigr)
\]
gives 
\((\sigma u,\vartheta\varphi)
=(\sigma\overline{\partial} u+\overline{\partial}\sigma\wedge u,\varphi)\), which is the
claimed distributional identity.
The right-hand side belongs to \(L^2\), since
\(\sigma\overline{\partial} u\in L^2\), \(\overline{\partial}\sigma\in L^\infty\), and
\(u\in L^2\). Hence \(\sigma u\in\dom(\overline{\partial})\) and the identity holds
in \(L^2\).

Now we prove that \(\sigma u\in\dom(\overline{\partial}^{\ast})\) by verifying the
Hilbert-adjoint pairing condition; by Lemma~\ref{lem:domains}(b) it
suffices to test against \(v\in\dom(\overline\partial_{(p,q-1)})\subset L^2_{(p,q-1)}(\Omega)\)
(for \(q=0\) there is nothing to prove). Let \(v\in\dom(\overline\partial_{(p,q-1)})\). Since
\(\sigma\in W^{1,\infty}(\Omega)\), the Leibniz identity gives
\(\sigma v\in\dom(\overline{\partial})\) with
\[
\overline{\partial}(\sigma v)
=
\overline{\partial}\sigma\wedge v+\sigma\,\overline{\partial} v .
\]
Using \(u\in\dom(\overline{\partial}^{\ast})\), we compute
\[
\begin{aligned}
(\sigma u,\overline{\partial} v)
&=
(u,\sigma\,\overline{\partial} v) \\
&=
\bigl(u,\overline{\partial}(\sigma v)-\overline{\partial}\sigma\wedge v\bigr) \\
&=
\bigl(\overline{\partial}^{\ast}u,\sigma v\bigr)
-
\bigl(u,\overline{\partial}\sigma\wedge v\bigr) \\
&=
\bigl(\sigma\,\overline{\partial}^{\ast}u,v\bigr)
-
\bigl(u,\overline{\partial}\sigma\wedge v\bigr).
\end{aligned}
\]
The last term is a continuous linear functional of
\(v\in L^2_{(p,q-1)}(\Omega)\), because
\(\overline{\partial}\sigma\in L^\infty_{\Lambda^{0,1}}\) and
\(u\in L^2_{(p,q)}\). Hence, by the Riesz representation theorem,
there exists
\(w_\mathrm{r}\in L^2_{(p,q-1)}(\Omega)\) such that
\[
(w_\mathrm{r},v)=\bigl(u,\overline{\partial}\sigma\wedge v\bigr)
\]
for all \(v\in L^2_{(p,q-1)}(\Omega)\). Consequently,
\[
(\sigma u,\overline{\partial} v)
=
\bigl(\sigma\,\overline{\partial}^{\ast}u-w_\mathrm{r},v\bigr)
\]
for all \(v\in\dom(\overline{\partial})\). By the definition of
\(\dom(\overline{\partial}^{\ast})\), this proves that
\[
\sigma u\in\dom(\overline{\partial}^{\ast}),
\qquad
\overline{\partial}^{\ast}(\sigma u)=\sigma\,\overline{\partial}^{\ast}u-w_\mathrm{r} .
\]

By the defining adjoint relation
\[
(u,\overline{\partial}\sigma\wedge v)
=
\bigl(\iota_{(\overline{\partial}\sigma)^\sharp}u,v\bigr)
\qquad
\text{for all }v,
\]
we identify
\[
w_\mathrm{r}=\iota_{(\overline{\partial}\sigma)^\sharp}u.
\]
so we have
\[
\overline{\partial}^{\ast}(\sigma u)
=
\sigma\,\overline{\partial}^{\ast}u-\iota_{(\overline{\partial}\sigma)^\sharp}u .
\]

By the domain-equality property,
\[
\dom(\overline{\partial}_f)=\dom(\overline{\partial}),
\qquad
\dom(\overline{\partial}_f^{\ast})=\dom(\overline{\partial}^{\ast})
\qquad\text{on every bidegree},
\]
we conclude
\[
\sigma u\in\dom(\overline{\partial}_f)\cap\dom(\overline{\partial}_f^{\ast}).
\]
Moreover, since \(\dom(\overline{\partial}_f^{\ast})=\dom(\overline{\partial}^{\ast})\) and the weak
\(\overline{\partial}\)-Neumann condition is automatic for elements of
\(\dom(\overline{\partial}^{\ast})\), the form \(\sigma u\) satisfies the weak
\(\overline{\partial}\)-Neumann condition, hence the weak \(\overline{\partial}_f\)-Neumann condition
by the equivalence.

Finally, using \(\overline{\partial}_f=\overline{\partial}+\partial f\wedge\) and
\(\overline{\partial}_f^{\ast}=\overline{\partial}^{\ast}+\iota_{(\partial f)^\sharp}\), we obtain
\[
\overline{\partial}_f(\sigma u)
=
\sigma \overline{\partial}_fu+\overline{\partial}\sigma\wedge u,
\qquad
\overline{\partial}_f^{\ast}(\sigma u)
=
\sigma \overline{\partial}_f^{\ast}u-\iota_{(\overline{\partial}\sigma)^\sharp}u
\]
in \(L^2\). This proves the product rules \eqref{eq:weak-product} and
completes the proof.
\end{proof}

\begin{lemma}[Relative localization in a domain]
\label{lem:relative-localization}
Let \(\Omega\subset\mathbb C^n\) be a bounded domain,
\(0<\epsilon<1\), and let \(U\subset\Omega\) be open.
Let \(g\in L^2(\Omega)\) be such that its restriction \(g|_U\) belongs to
\(H^\epsilon(U)\), and assume
$\operatorname{supp}g\subset U$,
$d:=\operatorname{dist}\bigl(\operatorname{supp}g,\Omega\setminus U\bigr)>0$.
Then \(g\in H^\epsilon(\Omega)\), and there exists a constant
$C=C(\epsilon,n,d)$
such that
\[
\|g\|_{H^\epsilon(\Omega)}
\le C\|g\|_{H^\epsilon(U)}.
\]
\end{lemma}

\begin{proof}
Since \(\operatorname{supp}g\subset U\) and
\(d=\operatorname{dist}(\operatorname{supp}g,\Omega\setminus U)>0\),
we have \(g=0\) a.e.\ on \(\Omega\setminus U\). Hence
$\|g\|_{L^2(\Omega)}=\|g\|_{L^2(U)}$.
Use the Gagliardo norm of Remark~\ref{rem:Heps-norms},
\[
\|g\|_{H^\epsilon(\Omega)}^2
=
\|g\|_{L^2(\Omega)}^2
+
\iint_{\Omega\times\Omega}
\frac{|g(x)-g(y)|^2}{|x-y|^{2n+2\epsilon}}\,dx\,dy .
\]
Because \(g\) vanishes a.e.\ on \(\Omega\setminus U\), the Gagliardo
integral receives contributions only from pairs \((x,y)\) with at
least one of \(x,y\) in \(U\).  By symmetry of the integrand, the two
cross regions \(U\times(\Omega\setminus U)\) and
\((\Omega\setminus U)\times U\) contribute equally, while the region
\((\Omega\setminus U)\times(\Omega\setminus U)\) contributes nothing.
Thus the integral equals the \(U\times U\) contribution plus twice the
cross contribution.
Then the first term is bounded by the Gagliardo
seminorm of \(g|_U\) on \(U\):

\[\iint_{U\times U}
\frac{|g(x)-g(y)|^2}{|x-y|^{2n+2\epsilon}}\,dx\,dy
=
\bigl[g|_U\bigr]_{H^\epsilon(U)}^2 .
\]
For the cross term, use \(g=0\) a.e.\ on \(\Omega\setminus U\) and
\(|x-y|\ge d\) for \(x\in\operatorname{supp}g\),
\(y\in\Omega\setminus U\). Thus
\[
\begin{aligned}
\iint_{U\times(\Omega\setminus U)}
\frac{|g(x)-g(y)|^2}{|x-y|^{2n+2\epsilon}}\,dx\,dy
&=
\iint_{\operatorname{supp}g\times(\Omega\setminus U)}
\frac{|g(x)|^2}{|x-y|^{2n+2\epsilon}}\,dx\,dy \\
&\le
\int_{\operatorname{supp}g}|g(x)|^2
\left(\int_{|x-y|\ge d}\frac{dy}{|x-y|^{2n+2\epsilon}}\right)dx \\
&\le
C(d,\epsilon,n)\|g\|_{L^2(U)}^2 .
\end{aligned}
\]
Combining these estimates with
\(\|g\|_{L^2(\Omega)}=\|g\|_{L^2(U)}\), we obtain

\[\|g\|_{H^\epsilon(\Omega)}^2
\le
C(\epsilon,n,d)
\bigl(
\|g\|_{L^2(U)}^2
+
\bigl[g|_U\bigr]_{H^\epsilon(U)}^2
\bigr)
\le
C(\epsilon,n,d)\|g\|_{H^\epsilon(U)}^2 .\]
\end{proof}

\begin{proposition}[Weighted subelliptic estimate on a boundary collar piece via Straube's local estimate]
\label{prop:subelliptic-collar-piece-L2}
Let \(\Omega\Subset\mathbb C^n\) be a bounded pseudoconvex domain with
Lipschitz boundary satisfying the boundary geometry hypothesis, with
uniform exponent \(\epsilon\in(0,1/2)\). Fix \(0\le p\le n\),
$1\le q\le n$.  Let \(V\subset\mathbb C^n\) be open with
\(V\cap\partial\Omega\ne\varnothing\) and \(\overline V\cap Z=\varnothing\),
and assume there are finitely many open sets \(U_1,\dots,U_K\) and
constants \(C_1,\dots,C_K>0\) such that
\(\overline V\subset\bigcup_{j=1}^KU_j\) and, for every \(j\) and every
\(\varphi\in L^2_{(0,q)}(\Omega)\cap\dom(\overline{\partial})\cap\dom(\overline{\partial}^{\ast})\),
\begin{equation}\label{eq:local-estimate-Uj-prop}
\|\varphi\|_{H^\epsilon(U_j\cap\Omega)}
\le C_j\Bigl(
\|\overline{\partial}\varphi\|_{L^2(\Omega)}
+\|\overline{\partial}^{\ast}\varphi\|_{L^2(\Omega)}
\Bigr).
\end{equation}
Let \(\theta\in W^{1,\infty}(\Omega)\) be the order-adapted weight and
$u\in L^2_{(p,q)}(\Omega)\cap \operatorname{dom}(\overline{\partial}_f)\cap
\operatorname{dom}(\overline{\partial}_f^{\ast})$.
Then
\(\theta u\in L^2_{(p,q)}(\Omega)\cap\operatorname{dom}(\overline{\partial}_f)\cap \operatorname{dom}(\overline{\partial}_f^{\ast})\),
and there exist an open set \(W\subset V\) with
\(V\cap\partial\Omega\subset W\), independent of \(u\), and a constant
\(C_V\), depending only on \(V\), \(W\), the sets \(U_j\), the constants
\(C_j\) and \(n,\epsilon\), such that
\begin{equation}
\label{eq:subelliptic-collar-piece-L2}
\|\theta u\|_{H^{\epsilon}(W\cap\Omega)}^2
\le C_V\Bigl(
\|\overline{\partial}_f(\theta u)\|_{L^2(\Omega)}^2
+\|\overline{\partial}_f^{\ast}(\theta u)\|_{L^2(\Omega)}^2
+\|\theta u\|_{L^2(\Omega)}^2
\Bigr).
\end{equation}
\end{proposition}

\begin{proof}
By Lemma~\ref{lem:multiplier-closure-L2} with \(\sigma=\theta\),
\(w:=\theta u\in L^2_{(p,q)}(\Omega)\cap\dom(\overline{\partial}_f)\cap\dom(\overline{\partial}_f^{\ast})\), and by
Lemma~\ref{lem:domains}(a),(b),
\(w\in\dom(\overline\partial)\cap\dom(\overline\partial^{\ast})\) with
\(\overline\partial w=\overline{\partial}_fw-\partial f\wedge w\) and
\(\overline\partial^{\ast}w=\overline{\partial}_f^{\ast}w-\iota_{(\partial f)^\sharp}w\).

\emph{The \((p,q)\) local estimate.}  By the consequence part of
Lemma~\ref{lem:componentwise}, \eqref{eq:local-estimate-Uj-prop} implies
\[
\|v\|_{H^\epsilon(U_j\cap\Omega)}
\le \sqrt2\,C_j\Bigl(\|\overline{\partial} v\|_{L^2(\Omega)}+\|\overline{\partial}^{\ast}v\|_{L^2(\Omega)}\Bigr)
\]
for every \(v\in L^2_{(p,q)}(\Omega)\cap\dom(\overline{\partial})\cap\dom(\overline{\partial}^{\ast})\),
with a constant independent of \(p\).

\emph{Choice of \(W\).}  Since \(\overline V\cap\partial\Omega\) is compact
and contained in the open set \(\bigcup_jU_j\),
\[
d_{\partial}:=\operatorname{dist}\Bigl(\overline{V}\cap\partial\Omega,\
\mathbb C^n\setminus\bigcup_jU_j\Bigr)>0 .
\]
Put
\(W:=V\cap\bigl\{z:\operatorname{dist}(z,\overline V\cap\partial\Omega)<d_\partial/2\bigr\}\).
Then \(W\) is open, independent of \(u\), contains \(V\cap\partial\Omega\),
and \(\overline W\subset\bigcup_jU_j\).  Choose
\(\chi_j\in C^\infty_c(U_j)\), \(0\le\chi_j\le1\), with
\(\sum_j\chi_j=1\) on a neighborhood of \(\overline W\).

\emph{Localization.}  By Lemma~\ref{lem:multiplier-closure-L2},
\(\chi_jw\in L^2_{(p,q)}(\Omega)\cap\dom(\overline{\partial})\cap\dom(\overline{\partial}^{\ast})\) with
\(\overline\partial(\chi_jw)=\chi_j\overline\partial w+\overline\partial\chi_j\wedge w\) and
\(\overline\partial^{\ast}(\chi_jw)=\chi_j\overline\partial^{\ast}w-\iota_{(\overline\partial\chi_j)^\sharp}w\).
Since \(\supp(\chi_jw)\subset\supp\chi_j\cap\overline\Omega\) and
\(\operatorname{dist}(\supp\chi_j,\mathbb C^n\setminus U_j)>0\),
Lemma~\ref{lem:relative-localization} (with \(U=U_j\cap\Omega\),
\(g=\chi_jw\)) gives
\(\|\chi_jw\|_{H^\epsilon(\Omega)}\le C^{\mathrm{loc}}_j\|\chi_jw\|_{H^\epsilon(U_j\cap\Omega)}\),
and a fortiori
\(\|\chi_jw\|_{H^\epsilon(W\cap\Omega)}\le C^{\mathrm{loc}}_j\|\chi_jw\|_{H^\epsilon(U_j\cap\Omega)}\).
As \(w=\sum_j\chi_jw\) on \(W\cap\Omega\), the triangle inequality for the
Gagliardo norm and the \((p,q)\) local estimate give
\begin{align*}
\|w\|_{H^\epsilon(W\cap\Omega)}^2
&\le K\sum_j(C^{\mathrm{loc}}_j)^2\|\chi_jw\|_{H^\epsilon(U_j\cap\Omega)}^2
\le 4K\sum_j(C^{\mathrm{loc}}_jC_j)^2
\bigl(\|\overline\partial(\chi_jw)\|^2+\|\overline\partial^{\ast}(\chi_jw)\|^2\bigr)\\
&\le C\Bigl(
\|\overline\partial w\|_{L^2(\Omega)}^2
+\|\overline\partial^{\ast} w\|_{L^2(\Omega)}^2
+\|w\|_{L^2(\Omega)}^2
\Bigr),
\end{align*}
where the \(L^2\)-term bounds the commutators
\(\overline\partial\chi_j\wedge w\) and \(\iota_{(\overline\partial\chi_j)^\sharp}w\) by
\(\|\nabla\chi_j\|_{L^\infty}\|w\|\).

\emph{Bridge.}  Since \(w\) has bidegree \((p,q)\), the four summands of
\(\overline{\partial}_fw\) and \(\overline{\partial}_f^{\ast}w\) lie in distinct bidegrees, so pointwise
\(|\overline\partial w|\le|\overline{\partial}_fw|\) and \(|\overline\partial^{\ast} w|\le|\overline{\partial}_f^{\ast}w|\).
Substituting and recalling \(w=\theta u\) gives
\eqref{eq:subelliptic-collar-piece-L2}.
\end{proof}

\section{Interior \texorpdfstring{\(H^{1/2}\)}{H 1/2} Estimate at an Isolated Critical Point}
\label{sec:critical}

Let \(z_c\in Z\), so \(m:=m_c\ge 1\), and fix \(0<R\le R_c\), so that
\(\theta=r_c^{m}\) on \(B_R(z_c)\).  (In the global assembly \(R=R_c\) and
the form under consideration is \(\chi_cu\), which is compactly supported
in \(B_R(z_c)\).)  We work in coordinates centered at \(z_c\):
\[
\tilde{z}:=z-z_c,\qquad r:=|\tilde{z}|,\qquad \varkappa:=\tilde{z}/r\quad (\tilde{z}\ne0),\qquad
\theta=r^m,\qquad \nu=2n,\qquad \beta=m+n-1,
\]
and
\[
\varsigma_R(r):=1+\log\frac{R}{r},\qquad 0<r\le R .
\]
In the sequel, to simplify the notation, we drop the tilde and write
\(z\) for the centered coordinate \(\tilde z\).

Throughout this section, a form \(u\in L^2_{(p,q)}(B_R(z_c))\) with
\(\supp u\Subset B_R(z_c)\) is identified with its zero extension to
\(\Omega\).  Accordingly, \(\dom(\overline\partial_f)\) and
\(\dom(\overline\partial_f^{\ast})\) always refer to the total-degree
realizations on \(\Omega\).  By Lemma~\ref{lem:domains}(c), for such a
compactly supported form this membership is equivalent to membership in
the corresponding domains formed on \(B_R(z_c)\); in particular
\(\overline\partial^{\ast}u=\vartheta u\) on \(B_R(z_c)\).  No ambiguity
arises from writing the domains without specifying the underlying open
set.

With \(\mathcal D=\overline\partial+\overline\partial^{\ast}\) the flat Dolbeault--Dirac
operator, locally put
\[
E_{r,\mathcal D}^0(u):=\int_{B_R(z_c)}r^{2m}|\mathcal D u|^2,
\,\,
E_r^f(u):=\int_{B_R(z_c)}r^{2m}\bigl(|\overline{\partial}_fu|^2+|\overline{\partial}_f^{\ast}u|^2\bigr),
\qquad
M_r(u):=\int_{B_R(z_c)}r^{2m}|u|^2 .
\]
For \(u\) of fixed bidegree \((p,q)\), pointwise orthogonality gives
\(|\mathcal Du|^2=|\overline\partial u|^2+|\overline\partial^{\ast}u|^2\le|\overline{\partial}_fu|^2+|\overline{\partial}_f^{\ast}u|^2\),
hence \(E_{r,\mathcal D}^0(u)\le E_r^f(u)\).

\begin{lemma}[Weighted gradient--Dirac multiplier inequality]
\label{lem:GD}
Let \(w\in H^1_0(B_R(z_c))\) be a compactly supported section of the Clifford
module, and let \(h:[0,R]\to[0,\infty)\) be a radial Lipschitz function.
Then there exists a constant \(C_{\mathrm{GD}}>0\), depending only on the
normalization of \(\mathcal D\), such that
\begin{equation}\label{GD}
\|h\nabla w\|_{L^2(B_R(z_c))}
\le
C_{\mathrm{GD}}\Bigl(
\|h\mathcal Dw\|_{L^2(B_R(z_c))}
+
\|h'w\|_{L^2(B_R(z_c))}
\Bigr).
\end{equation}
\end{lemma}

\begin{proof}
For compactly supported \(H^1\) sections \(g\) of the full form bundle,
integration by parts with the flat identity
\((\overline\partial+\overline\partial^{\ast})^2=-\frac14\Delta\) of the
inner-product convention gives
$\|\nabla g\|_2=2\,\|\mathcal Dg\|_2 $.
Since \(h\) is radial and Lipschitz, \(\nabla h=h'(r)\,\varkappa\) a.e., and the
principal symbol of \(\mathcal D\) gives the product rule
\[
\mathcal D(hw)=h\mathcal Dw+\frac{1}{2}h'\,U(\varkappa)w,
\]
where \(U(\varkappa):=\sum_{a=1}^{\nu}c_a\varkappa_a\) is Clifford multiplication by the
radial unit vector \(\varkappa=x/r\), which is unitary because
\(U(\varkappa)^{\ast}=-U(\varkappa)\) and \(U(\varkappa)^2=-|\varkappa|^2=-1\); hence
\[
\bigl|[\mathcal D,h]w\bigr|
=\frac{1}{2}|h'|\,|w|
\le |h'|\,|w| .
\]
Consequently
\[
\|h\nabla w\|_2
\le\|\nabla(hw)\|_2+\|h'w\|_2
=2\|\mathcal D(hw)\|_2+\|h'w\|_2
\le 2\|h\mathcal Dw\|_2+2\|h'w\|_2
\]
Thus \eqref{GD} holds with \(C_{\mathrm{GD}}:=2\).
The estimate for \(w\) in the stated class follows
by smooth approximation in \(H^1_0\).
\end{proof}

\begin{lemma}[\(H^{1/2}\) estimate via weighted gradient energy and Stein--Weiss]
\label{lem:SW}
Let $\nu\ge 2$. There is a constant $C_\nu$ such that for every compactly
supported $g\in H^1(\R^\nu)$ (scalar valued; for Clifford-module valued $g$
the inequality is applied componentwise and summed),
\begin{equation}
\label{eq:Dfs}
\|g\|_{H^{1/2}(\R^\nu)}^2
\le
C_\nu\Bigl(\|g\|_{L^2(\R^\nu)}^2+\int_{\R^\nu}|x|\,|\nabla g|^2\,dx\Bigr).
\end{equation}
\end{lemma}

\begin{proof}
We first assume \(g\in C^\infty_c(\R^\nu)\); then
\(\mathscr{G}:=(-\Delta)^{1/2}g\) is smooth with \(|\mathscr{G}(x)|=O(|x|^{-\nu-1})\) as
\(|x|\to\infty\),  here $\Delta:=\sum_{\ell=1}^\nu\partial_{x_\ell}^2$ is real Euclidean Laplacian. So all Fourier-multiplier and Riesz-potential
representations below are classical and agree.  For a general compactly
supported \(g\in H^1(\R^\nu)\), 
let
\(\{\eta_j\}_{j\ge1}\) be a standard mollifier and set
$g_j:=\eta_j\ast g$.
Then \(g_j\in C_c^\infty(\mathbb R^\nu)\), all supports \(\supp g_j\) lie
in a fixed ball (independent of \(j\)), and
$g_j\to g$ in $H^1(\mathbb R^\nu)$, hence
\(\|g_j\|_{H^{1/2}}\to\|g\|_{H^{1/2}}\), \(\|g_j\|_{L^2}\to\|g\|_{L^2}\),
and \(\int|x||\nabla g_j|^2\to\int|x||\nabla g|^2\) (the weight \(|x|\)
being bounded on the common support); so \eqref{eq:Dfs} passes to the
limit.  It suffices to prove the homogeneous bound
$\|g\|_{\dot H^{1/2}}^2\le C_\nu\int_{\R^\nu}|x|\,|\nabla g|^2\,dx$, since
$\|g\|_{H^{1/2}}^2\simeq\|g\|_{L^2}^2+\|g\|_{\dot H^{1/2}}^2$. Because $g$ is
compactly supported and lies in $H^1(\R^\nu)$, we have
$\nabla g\in L^2$, $(-\Delta)^{1/2}g\in L^2$, and
$\|g\|_{\dot H^{1/2}}^2=\|(-\Delta)^{1/4}g\|_{L^2}^2$.

We use the identity
\(\bigl\||x|^{1/2}G\bigr\|_{L^2(\R^\nu)}^2=\|G\|_{L^2(\varpi\,dx)}^2\) with
\(\varpi(x)=|x|\); thus the weight \(|x|^{1/2}\) multiplying the function and the weight \(|x|\) in the measure are, when squared, the same statement. The Riesz
potential \(I_{1/2}\) below is viewed as the map
\(I_{1/2}:L^2(|x|\,dx)\to L^2(dx)\), and the Riesz transforms \(\mathcal{R}_\ell\) as
bounded maps \(L^2(|x|\,dx)\to L^2(|x|\,dx)\).

\emph{Stein--Weiss.}
Factor the Fourier multiplier $|\xi|^{1/2}=|\xi|^{-1/2}\cdot|\xi|$, i.e.
$(-\Delta)^{1/4}g=I_{1/2}\,\mathscr{G}$,
where $I_{\gamma}$ is the Riesz potential of order $\gamma$, with Fourier
multiplier $|\xi|^{-\gamma}$. The Stein--Weiss inequality
\cite{SteinWeiss} states that for $0<\gamma<\nu$ and
$1<\mathfrak{p}\le \mathfrak{q}<\infty$,
\[
\bigl\||x|^{-a}I_\gamma \mathscr{G}\bigr\|_{L^{\mathfrak{q}}(\R^\nu)}
\le
 C_{\mathrm{SW}}\,\bigl\||x|^{b}\mathscr{G}\bigr\|_{L^{\mathfrak{p}}(\R^\nu)}
\]
provided $a+b\ge0$, $a<\nu/\mathfrak{q}$, $b<\nu/\mathfrak{p}^{\ast}$,  and the
dimensional balance
$\frac{1}{\mathfrak{p}}-\frac{1}{\mathfrak{q}}=\frac{\gamma-a-b}{\nu}$ holds, where $\mathfrak{p}^{\ast}=\mathfrak{p}/(\mathfrak{p}-1)$ is H\"{o}lder conjugate of $\mathfrak{p}$. Take
\[
\mathfrak{p}=\mathfrak{q}=2,\qquad \gamma=\tfrac12,\qquad a=0,\qquad b=\tfrac12 .
\]
Then $\gamma=a+b$ gives the balance $\frac1{\mathfrak{p}}-\frac1{\mathfrak{q}}=0$;
$a+b=\tfrac12\ge0$; $a=0<\nu/\mathfrak{q}=\nu/2$; and
$b=\tfrac12<\nu/\mathfrak{p}^{\ast}=\nu/2$ because $\nu\ge2$.

In the original formulation \cite[Theorem~\(B^\ast\)]{SteinWeiss}, with
\(T_{\lambda_{\rm SW}} \mathscr{G}(x)=\int|x-y|^{-\lambda_{\rm SW}}\mathscr{G}(y)\,dy\), one has
\(I_\gamma=c_{\nu,\gamma}T_{\nu-\gamma}\); the above choice corresponds to
\[
\lambda_{\rm SW}=\nu-\tfrac12\in(0,\nu),\qquad
\mathfrak p=\mathfrak q=2,
\]
input weight exponent
$\alpha_{\rm in}=b=\tfrac12<\nu/\mathfrak p^{\ast}$,
output weight exponent
$\beta_{\rm out}=a=0<\nu/\mathfrak q$,
with
\[
\alpha_{\rm in}+\beta_{\rm out}\ge0,
\qquad
\frac1{\mathfrak q}
=
\frac1{\mathfrak p}
+
\frac{\lambda_{\rm SW}+\alpha_{\rm in}+\beta_{\rm out}}{\nu}
-1,
\]
since
$\lambda_{\rm SW}+\alpha_{\rm in}+\beta_{\rm out}=\nu$.
The right-hand side is finite, because
$\mathscr{G}$ is smooth and $|\mathscr{G}(x)|=O(|x|^{-\nu-1})$, so that $|x|^{1/2}\mathscr{G}\in L^2$. Hence
\begin{equation}
\label{eq:SWapp}
\|(-\Delta)^{1/4}g\|_{L^2}
=\|I_{1/2}\mathscr{G}\|_{L^2}
\le C_{\mathrm{SW}}(\nu)\bigl\||x|^{1/2}\mathscr{G}\bigr\|_{L^2}
=C_{\mathrm{SW}}(\nu)\bigl\||x|^{1/2}(-\Delta)^{1/2}g\bigr\|_{L^2}.
\end{equation}

\emph{Weighted Riesz bound.}
Write $(-\Delta)^{1/2}=\sum_{\ell=1}^{\nu}\mathcal{R}_\ell\,\partial_{x_\ell}$, where
$\mathcal{R}_\ell=-\partial_{x_\ell}(-\Delta)^{-1/2}$ are the Riesz transforms (Fourier
multipliers $-i\xi_\ell/|\xi|$); indeed
$\sum_\ell(-i\xi_\ell/|\xi|)(i\xi_\ell)=|\xi|$. 
The weight \(\varpi(x)=|x|\) is a
Muckenhoupt \(A_2\) weight on \(\R^\nu\), since \(|x|^{\mathfrak{m}}\in A_2(\R^\nu)\) exactly
when \(-\nu<\mathfrak{m}<\nu\) \cite[Example 7.1.7]{Grafakos2014}, and here \(\mathfrak{m}=1\in(-\nu,\nu)\) for \(\nu\ge2\); its \(A_2\)
characteristic \([\varpi]_{A_2}\) is finite and depends only on \(\nu\).
\footnote{Here and below, for a weight \(w\) and a cube \(Q\), we use the
standard Muckenhoupt characteristics with cube averages:
\[
[w]_{A_p}
:=
\sup_{Q}
\Bigl(
\frac{1}{|Q|}\int_Q w\,dx
\Bigr)
\Bigl(
\frac{1}{|Q|}\int_Q w^{-1/(p-1)}\,dx
\Bigr)^{p-1},
\qquad
[w]_{A_2}
:=
\sup_{Q}
\Bigl(
\frac{1}{|Q|}\int_Q w\,dx
\Bigr)
\Bigl(
\frac{1}{|Q|}\int_Q w^{-1}\,dx
\Bigr).
\]
The weighted bound in \cite[Theorem~1.3]{Hyt} is stated for this
characteristic.}
Each $\mathcal{R}_\ell$
is a Calder\'on--Zygmund operator, so by the weighted bound for such
operators \cite[Theorem~1.3]{Hyt}, applied on \(\R^\nu\) to the
Calder\'on--Zygmund kernel $c_\nu(x_\ell-y_\ell)|x-y|^{-\nu-1}$ (H\"older
exponent $\alpha=1$) and to $F\in C^\infty_c\subset L^2(\R^\nu)\cap L^2(\varpi)$
(boundedness here is classical, due to
Hunt--Muckenhoupt--Wheeden \cite{HuMuckeWhe} and Coifman--Fefferman
\cite{CoifFeffer}),
\[
\bigl\||x|^{1/2}\mathcal{R}_\ell F\bigr\|_{L^2}
=\|\mathcal{R}_\ell F\|_{L^2(\varpi)}
\le C_{\mathrm{CZ}}(\nu)\,[\varpi]_{A_2}\,\|F\|_{L^2(\varpi)}
=C_{\mathrm{Riesz}}(\nu)\,\bigl\||x|^{1/2}F\bigr\|_{L^2}.
\]
Applying this with $F=\partial_{x_\ell} g$ and summing over $\ell$,
\begin{equation}
\label{eq:Rieszapp}
\bigl\||x|^{1/2}(-\Delta)^{1/2}g\bigr\|_{L^2}
\le\sum_{\ell}\bigl\||x|^{1/2}\mathcal{R}_\ell(\partial_{x_\ell} g)\bigr\|_{L^2}
\le C_{\mathrm{Riesz}}(\nu)\,\bigl\||x|^{1/2}\nabla g\bigr\|_{L^2}.
\end{equation}

Combining \eqref{eq:SWapp} and \eqref{eq:Rieszapp},
\[
\|g\|_{\dot H^{1/2}}^2
=\|(-\Delta)^{1/4}g\|_{L^2}^2
\le \nu\,C_{\mathrm{SW}}(\nu)^2C_{\mathrm{Riesz}}(\nu)^2\,\int_{\R^\nu}|x|\,|\nabla g|^2\,dx,
\]
and adding $\|g\|_{L^2}^2$ yields \eqref{eq:Dfs}. For Clifford-module
valued $g$, the operators act componentwise and the bound is summed over
components.
\end{proof}

\begin{lemma}[Angular spectrum and spectral gap for the Dolbeault--Dirac module]
\label{lem:angular-spectrum-gap}
Let \(\nu=2n\), and let
\(\mathcal D=\overline{\partial}+\overline{\partial}^{\ast}\) act on the full form bundle
\(\Lambda^\bullet\mathbb C^n\). In polar coordinates \(z=r\varkappa\),
\(r>0\), \(\varkappa\in S^{\nu-1}\), there is a self-adjoint
first-order operator \(\mathcal L\) on \(L^2(S^{\nu-1},\Lambda^\bullet)\)
such that
\[
\mathcal D=\frac{1}{2}U(\varkappa)\Bigl(\partial_r+\frac1r \mathcal L\Bigr),
\]
where \(U(\varkappa)\) is unitary Clifford multiplication by the radial
vector. Then the Clifford identities
\[
\mathcal L^2-(\nu-2)\mathcal L=-\Delta_{S^{\nu-1}}I,
\qquad
\mathcal LU(\varkappa)+U(\varkappa)\mathcal L=(\nu-1)U(\varkappa)
\]
hold. Moreover, the spectrum of \(\mathcal L\) consists of integers:
$\operatorname{spec}(\mathcal L)\subset\mathbb Z$.
Let \(m\ge1\) and
\[
\beta=m+n-1=\nu-1+(m-n).
\]
Then \(\beta\in\mathbb Z_{\ge1}\), and for every eigenvalue
\(\mu\ne\beta\) of \(\mathcal L\),
$|\mu-\beta|\ge1$.
\end{lemma}

\begin{proof}
The polar decomposition and the Clifford identities are standard for the
Euclidean Dirac operator on the full exterior algebra (cf.\ \cite[Chapter~II]{LawsonMichelsohn}).  We recall the
short derivation.  By the inner-product convention, put
\[
\mathfrak D:=2\,\mathcal D=\sum_{a=1}^{\nu}c_a\partial_{x_a},
\qquad
\mathcal L:=-\sum_{a<b}c_ac_b(x_a\partial_b-x_b\partial_a),
\]
where the constant Clifford generators \(c_a\) are skew-adjoint and
satisfy \(c_ac_b+c_bc_a=-2\delta_{ab}\), so that \(\mathfrak D^2=-\Delta\).
With \(x=\sum_a x_ac_a\), \(r=|x|\), \(U(\varkappa):=x/r\) (which is unitary,
since \(U(\varkappa)^{\ast}=-U(\varkappa)\) and \(U(\varkappa)^2=-1\)), a direct
computation using \(c_a^2=-1\) gives
\[
x\,\mathfrak D=-r\partial_r-\mathcal L,
\qquad
\mathcal Lx+x\mathcal L=(\nu-1)x .
\]
Multiplying the first identity on the left by \(x/r^2\) and using
\(x^2=-r^2\) gives
\[
\mathfrak D=U(\varkappa)\Bigl(\partial_r+\frac1r\mathcal L\Bigr),
\qquad\text{hence}\qquad
\mathcal D=\frac12U(\varkappa)\Bigl(\partial_r+\frac1r\mathcal L\Bigr).
\]
For the second identity, let \(\mathfrak{s}\) be a smooth section.  Since
\(X_{ab}x_d=x_a\delta_{bd}-x_b\delta_{ad}\),
\[
\mathcal L(x\mathfrak{s})
=-\sum_{a<b}c_ac_b(x_ac_b-x_bc_a)\mathfrak{s}-\sum_{a<b}c_ac_b\,x\,X_{ab}\,\mathfrak{s}
=(\nu-1)x\mathfrak{s}-\sum_{a<b}c_ac_b\,x\,X_{ab}\,\mathfrak{s},
\]
using \(c_ac_bc_b=-c_a\) and \(c_ac_bc_a=c_b\) for \(a\ne b\).  Hence
\(\mathcal Lx+x\mathcal L=(\nu-1)x-\sum_{a<b}(c_ac_bx+xc_ac_b)X_{ab}\).
For \(d\in\{a,b\}\) the generator \(c_d\) anticommutes with \(c_ac_b\), and
for \(d\notin\{a,b\}\) it commutes, so
\(\sum_{a<b}(c_ac_bx+xc_ac_b)X_{ab}
=\sum_{(d,a,b)\ \mathrm{distinct}}c_dc_ac_b\,x_dx_a\partial_b\),
which vanishes because \(c_dc_ac_b\) is antisymmetric and \(x_dx_a\)
symmetric in \((d,a)\).  This proves \(\mathcal Lx+x\mathcal L=(\nu-1)x\).
Since \(\mathcal L\) contains no radial derivative, it commutes with
multiplication by \(r\), and dividing by \(r\) gives
\(\mathcal LU(\varkappa)+U(\varkappa)\mathcal L=(\nu-1)U(\varkappa)\).  Consequently
\(\bigl(\partial_r+\frac1r\mathcal L\bigr)U(\varkappa)
=U(\varkappa)\bigl(\partial_r+\frac{\nu-1-\mathcal L}{r}\bigr)\), and
\[
-\Delta=\mathfrak D^2
=-\Bigl(\partial_r+\frac{\nu-1-\mathcal L}{r}\Bigr)
\Bigl(\partial_r+\frac{\mathcal L}{r}\Bigr)
=-\Bigl(\partial_r^2+\frac{\nu-1}{r}\partial_r
+\frac{(\nu-2)\mathcal L-\mathcal L^2}{r^2}\Bigr).
\]
Comparing with
\(\Delta=\partial_r^2+\frac{\nu-1}{r}\partial_r+\frac1{r^2}\Delta_{S^{\nu-1}}\)
(acting componentwise in the constant frame) yields
\[
\mathcal L^2-(\nu-2)\mathcal L=-\Delta_{S^{\nu-1}}I .
\]

\emph{Symmetry of \(\mathcal L\).}
The vector fields
\[
X_{ab}:=x_a\partial_b-x_b\partial_a,\qquad 1\le a<b\le\nu,
\]
are Killing fields of \(S^{\nu-1}\); being divergence free on the
sphere, they are skew-symmetric on \(C^\infty(S^{\nu-1},\Lambda^\bullet)\)
with respect to the \(L^2(S^{\nu-1})\) inner product.  For \(a\ne b\)
the constant matrix \(c_ac_b\) satisfies
\((c_ac_b)^{\ast}=c_b^{\ast}c_a^{\ast}=c_bc_a=-c_ac_b\), and it commutes with \(X_{ab}\).
Hence each summand \(c_ac_bX_{ab}\) is symmetric, and so is
\(\mathcal L\) on \(C^\infty(S^{\nu-1},\Lambda^\bullet)\).
(The self-adjointness of \(\mathcal L\) is not inferred from that of
\(\mathcal D\): in the polar formula the factor \(U(\varkappa)\) is a
variable unitary multiplier and \(\partial_r^{\ast}=-\partial_r-\frac{\nu-1}{r}\)
in the measure \(r^{\nu-1}dr\,d\varkappa\).)

\emph{Spectral decomposition.}
For \(k\ge0\) let \(H_k\) be the (finite-dimensional) space of
restrictions to \(S^{\nu-1}\) of \(\Lambda^\bullet\mathbb C^n\)-valued
harmonic polynomials homogeneous of degree \(k\).  The spaces \(H_k\) are
mutually orthogonal, their algebraic sum is dense in
\(L^2(S^{\nu-1},\Lambda^\bullet)\), and
\(-\Delta_{S^{\nu-1}}=k(k+\nu-2)\) on \(H_k\) \cite[Chapter 3]{ArfReBM}.  Each \(X_{ab}\)
commutes with the Euclidean Laplacian and preserves homogeneity, and the
\(c_ac_b\) are constant; hence \(\mathcal L(H_k)\subset H_k\).  Since
\(\mathcal L\) is symmetric on the finite-dimensional space \(H_k\), it is
diagonalizable there with real eigenvalues and an orthonormal eigenbasis.
Collecting these bases over \(k\) gives an orthonormal basis of
\(L^2(S^{\nu-1},\Lambda^\bullet)\) consisting of smooth eigensections
of \(\mathcal L\); the closure of \(\mathcal L\) from
\(\bigoplus_k^{\mathrm{alg}}H_k\) is therefore self-adjoint, and for every
\(\mu\in\mathbb R\) the spectral projection \(\mathrm{P}_\mu\) is the orthogonal
projection onto the span of the basis elements with eigenvalue \(\mu\).
On \(H_k\) the first Clifford identity gives
\[
(\mathcal L+k)(\mathcal L-(\nu-2+k))=0,
\]
so every eigenvalue of \(\mathcal L|_{H_k}\) lies in
\(\{-k,\nu-2+k\}\).  On \(H_0\) (constants) one has
\(\mathcal L=0\), because \(\mathcal L\) is a first-order operator without
zero-order term.  Hence
\[
\operatorname{spec}(\mathcal L)
\subset
\{0,-1,-2,\dots\}\cup\{\nu-1,\nu,\dots\}
\subset\mathbb Z.
\]
Moreover, for \(k\ge1\) the two values \(-k\) and \(\nu-2+k\) are
distinct, so a positive eigenvalue \(\mu\) can occur only in the single
space \(H_{\mu-\nu+2}\); in particular every spectral projection
\(\mathrm{P}_\mu\) with \(\mu>0\) has finite rank and range consisting of smooth
sections.

Since \(m\ge1\), the number \(\beta=m+n-1\) is a positive integer.
Therefore \(\beta\) and any eigenvalue \(\mu\ne\beta\) of \(\mathcal L\)
differ by a nonzero integer, which gives \(|\mu-\beta|\ge1\).
\end{proof}

\begin{proposition}[Local \(H^{1/2}\) estimate at an isolated critical point: smooth case]
\label{prop:M-critical}
Let \(z_c\in Z\) with \(m:=m_c\ge1\). Let
\[
u\in C^{\infty}_{c,\Lambda^{p,q}}(B_R(z_c)),
\qquad
\operatorname{supp}u\Subset B_R(z_c),
\]
be a smooth, compactly supported \((p,q)\)-form.
Then there exists a constant
\(C=C(n,m,R)\), such that
\begin{equation}
\label{eq:M}
\|r^m u\|_{H^{1/2}(B_R(z_c))}^2
\le
C\Bigl(
\|r^m \overline{\partial}_fu\|_{L^2(B_R(z_c))}^2
+
\|r^m \overline{\partial}_f^{\ast}u\|_{L^2(B_R(z_c))}^2
+
\|r^m u\|_{L^2(B_R(z_c))}^2
\Bigr).
\end{equation}
Equivalently, in the notation above,
\[
\|r^m u\|_{H^{1/2}(B_R(z_c))}^2
\le
C\bigl(E_r^f(u)+M_r(u)\bigr).
\]

\end{proposition}

\begin{proof}

\medskip

\emph{Angular spectrum and the resonant gap.}
By Lemma~\ref{lem:angular-spectrum-gap},
\(\operatorname{spec}(\mathcal L)\subset\mathbb Z\), and for
\[
\beta=m+n-1
\]
every eigenvalue \(\mu\ne\beta\) of \(\mathcal L\) satisfies
\(|\mu-\beta|\ge1\). Let \(\mathrm{P}_\beta\) be the self-adjoint spectral
projection of \(\mathcal L\) onto the eigenvalue \(\beta\), with the
convention \(\mathrm{P}_\beta=0\) if
\(\beta\notin\operatorname{spec}(\mathcal L)\). Put
\[
u_\beta:=\mathrm{P}_\beta u,\qquad u_\perp:=(I-\mathrm{P}_\beta)u .
\]
Then either \(\beta\) is an eigenvalue and \(\mathrm{P}_\beta\ne0\), or
\(\beta\notin\operatorname{spec}(\mathcal L)\) and \(\mathrm{P}_\beta=0\); in
both cases the non-resonant channel satisfies the gap condition, and the
resonant channel is either non-trivial or empty.

Since \(u\in C_c^\infty(B_R(z_c))\) and \(\operatorname{supp}u\Subset B_R(z_c)\), there exists \(R'<R\) with \(\operatorname{supp}u\subset\overline{B_{R'}}\subset B_R(z_c)\), and \(u\in H^1_0(B_R(z_c))\). The spectral projection \(\mathrm{P}_\beta\) acts only on the angular variable and leaves the radial variable unchanged, so it does not enlarge the radial support; hence the radial projections of \(\operatorname{supp}u_\beta\) and \(\operatorname{supp}u_\perp\) are contained in the radial projection of \(\operatorname{supp}u\), and in particular
 \[
 \operatorname{supp}u_\beta,\ \operatorname{supp}u_\perp
 \subset \overline{B_{R'}}\subset B_R(z_c).
 \]
To conclude \(u_\beta,u_\perp\in H^1_0(B_R(z_c))\), it therefore remains only to show that they belong to \(H^1(B_R(z_c))\).

By Lemma~\ref{lem:angular-spectrum-gap}, \(\mathrm{P}_\beta\) has finite rank and
its range consists of smooth sections \(Y_1,\dots,Y_{K_\beta}\) (orthonormal);
hence \(u_\beta(r\varkappa)=\sum_{i=1}^{K_\beta}\varrho_i(r)Y_i(\varkappa)\) with
\(\varrho_i(r)=\langle u(r\,\cdot),Y_i\rangle_{L^2(S^{\nu-1})}\)
smooth on \([0,R)\), so \(u_\beta\) and \(u_\perp=u-u_\beta\) are smooth
on \(B_R(z_c)\setminus\{z_c\}\) and bounded.
We use the following removability fact: a form that is smooth and bounded on \(B_R(z_c)\setminus\{z_c\}\) and whose classical gradient is square integrable on this punctured ball belongs to \(H^1(B_R(z_c))\), and its weak gradient coincides almost everywhere on \(B_R(z_c)\) with the classical gradient.
Indeed, in \(\mathbb R^\nu\) with \(\nu\ge2\),
the singleton \(\{z_c\}\) is closed in \(B_R(z_c)\) and has zero \((\nu-1)\)-dimensional Hausdorff measure; therefore
it is removable for \(V_2^1=H^1\) by
\cite[Theorem~1.1.18]{Mazya2011}. Then it suffices to bound
the classical gradients on \(B_R(z_c)\setminus\{z_c\}\).

In polar coordinates, for every smooth section \(v\) one has
\[
|\nabla v|^2=|\partial_r v|^2+r^{-2}|\nabla_{S^{\nu-1}}v|^2 .
\]
Because \(\mathrm{P}_\beta\) commutes with \(\partial_r\), for almost every \(r\),
\[
\partial_r u_\beta=\mathrm{P}_\beta\partial_r u,
\qquad
\partial_r u_\perp=(I-\mathrm{P}_\beta)\partial_r u .
\]
Since \(\mathrm{P}_\beta\) and \(I-\mathrm{P}_\beta\) are orthogonal projections on
\(L^2(S^{\nu-1})\), they are \(L^2\)-contractions, and therefore
\[
\int_{S^{\nu-1}}|\partial_r u_\beta|^2\,d\varkappa
\le
\int_{S^{\nu-1}}|\partial_r u|^2\,d\varkappa,
\qquad
\int_{S^{\nu-1}}|\partial_r u_\perp|^2\,d\varkappa
\le
\int_{S^{\nu-1}}|\partial_r u|^2\,d\varkappa .
\]
For the angular part, recall that
\[
-\Delta_{S^{\nu-1}}I
=
\mathcal L^2-(\nu-2)\mathcal L .
\]
Thus \(-\Delta_{S^{\nu-1}}\) commutes with \(\mathrm{P}_\beta\), and
\(\mathrm{P}_\beta\) is the spectral projection onto the eigenspace
\(\mathcal L=\beta\).  Since \(\mathrm{P}_\beta\) has finite rank onto an eigenspace of \(-\Delta_{S^{\nu-1}}\) with eigenvalue \(\beta^2-(\nu-2)\beta\ge0\) (cf.\ the functional calculus \cite[Theorems~VIII.4--VIII.6]{ReedSimon}), \(\mathrm{P}_\beta\) commutes with
\((-\Delta_{S^{\nu-1}})^{1/2}\) and is an \(L^2\)-contraction on the
angular \(H^1\)-norm:
\[
\|\mathrm{P}_\beta v\|_{H^1(S^{\nu-1})}
\le
\|v\|_{H^1(S^{\nu-1})},
\qquad
\|(I-\mathrm{P}_\beta)v\|_{H^1(S^{\nu-1})}
\le
\|v\|_{H^1(S^{\nu-1})}.
\]
Equivalently,
\[
\int_{S^{\nu-1}}|\nabla_{S^{\nu-1}} u_\beta|^2\,d\varkappa
\le
\int_{S^{\nu-1}}|\nabla_{S^{\nu-1}} u|^2\,d\varkappa,
\]
and the same inequality holds with \(u_\beta\) replaced by \(u_\perp\).
Therefore
\[
\int_{B_R(z_c)}|\nabla u_\beta|^2\,dV
\le
\int_{B_R(z_c)}|\nabla u|^2\,dV<\infty,
\]
and similarly
\[
\int_{B_R(z_c)}|\nabla u_\perp|^2\,dV
\le
\int_{B_R(z_c)}|\nabla u|^2\,dV<\infty .
\]
Together with \(u_\beta,u_\perp\in L^2(B_R(z_c))\) and their compact
support in \(B_R(z_c)\), this yields
$u_\beta,\ u_\perp\in H^1_0(B_R(z_c))$.
In particular, no integrability of \(r^{-2}|u_\perp|^2\) is used here;
that quantity may fail to be integrable when \(n=1\) and \(u_\perp\)
contains a nonzero constant angular mode.

With this at hand, the cross terms in the gradient energy are well defined
and vanish by the orthogonality of \(\mathrm{P}_\beta\).

\emph{Orthogonal splitting of the gradient energy.} In a fixed Cartesian
trivialization of the Clifford module the gradient energy splits into radial
and angular parts,
\begin{equation}
\label{eq:polar-grad}
|\nabla u|^2=|\partial_r u|^2+\tfrac1{r^2}|\nabla_{S^{\nu-1}}u|^2,
\qquad
\int_{S^{\nu-1}}|\nabla_{S^{\nu-1}}u|^2\,d\varkappa
=\bigl\langle-\Delta_{S^{\nu-1}}u,u\bigr\rangle_{L^2(S^{\nu-1})} .
\end{equation}
By the first Clifford identity, $-\Delta_{S^{\nu-1}}I=\mathcal L^2-(\nu-2)\mathcal L$
is a polynomial in $\mathcal L$; therefore it commutes with the spectral
projection $\mathrm{P}_\beta$ and acts as the scalar
$\beta^2-(\nu-2)\beta$ on $\operatorname{ran}\mathrm{P}_\beta$. Consequently, at each fixed
radius $r$ the two cross terms vanish: 
since \(\mathrm{P}_\beta\) is a self-adjoint projection in \(L^2(S^{\nu-1})\) and commutes with \(\partial_r\), we have
 \[
 \partial_r u_\beta=\mathrm{P}_\beta\partial_r u,\qquad
 \partial_r u_\perp=(I-\mathrm{P}_\beta)\partial_r u,
 \]
 and hence
 \[
 \langle\partial_r u_\beta,\partial_r u_\perp\rangle
 =
 \langle \mathrm{P}_\beta\partial_r u,(I-\mathrm{P}_\beta)\partial_r u\rangle
 =0.
 \]
Likewise,
\[
\bigl\langle-\Delta_{S^{\nu-1}}u_\beta,u_\perp\bigr\rangle_{L^2(S^{\nu-1})}
=\bigl(\beta^2-(\nu-2)\beta\bigr)
\bigl\langle u_\beta,u_\perp\bigr\rangle_{L^2(S^{\nu-1})}=0 .
\]
Integrating \eqref{eq:polar-grad} against the radial density
$\kappa(r)\,r^{\nu-1}\,dr$ and using these two identities gives, for every
measurable radial weight $\kappa(r)\ge0$,
\begin{equation}
\label{eq:grad-split}
\int \kappa\,|\nabla u|^2\,dV=\int \kappa\,|\nabla u_\beta|^2\,dV+\int \kappa\,|\nabla u_\perp|^2\,dV,
\,\,
\int \kappa\,|u|^2\,dV=\int \kappa\,|u_\beta|^2\,dV+\int \kappa\,|u_\perp|^2\,dV .
\end{equation}
The second identity follows from the \(L^2(S^{\nu-1})\)-orthogonality
of \(u_\beta(r\,\cdot)\) and \(u_\perp(r\,\cdot)\) on each sphere.

Since \(u_\beta,u_\perp\in H^1_0(B_R(z_c))\) and the cross terms vanish,
the identities \eqref{eq:grad-split} hold for every radial weight
\(\kappa\ge0\) (both sides possibly \(+\infty\)). In particular, Lemma~\ref{lem:GD} applies to each channel.

\emph{Orthogonal splitting of the Dirac energy.}
Since \(U(\varkappa)\) is unitary, the polar decomposition of
Lemma~\ref{lem:angular-spectrum-gap} gives, pointwise on
\(B_R(z_c)\setminus\{z_c\}\),
\[
|\mathcal Dv|^2
=
\frac14\Bigl|\Bigl(\partial_r+\frac1r\mathcal L\Bigr)v\Bigr|^2
\qquad(v\in\{u,u_\beta,u_\perp\}).
\]
Because \(\mathrm{P}_\beta\) commutes with \(\partial_r\) and with \(\mathcal L\),
\[
\Bigl(\partial_r+\frac1r\mathcal L\Bigr)u_\beta
=
\mathrm{P}_\beta\Bigl(\partial_r+\frac1r\mathcal L\Bigr)u,
\qquad
\Bigl(\partial_r+\frac1r\mathcal L\Bigr)u_\perp
=
(I-\mathrm{P}_\beta)\Bigl(\partial_r+\frac1r\mathcal L\Bigr)u ,
\]
and since \(\mathrm{P}_\beta\) is an orthogonal projection on
\(L^2(S^{\nu-1})\), for almost every \(r\)
\[
\int_{S^{\nu-1}}|\mathcal Du|^2\,d\varkappa
=
\int_{S^{\nu-1}}|\mathcal Du_\beta|^2\,d\varkappa
+
\int_{S^{\nu-1}}|\mathcal Du_\perp|^2\,d\varkappa .
\]
Integrating in \(r\) with the radial weight \(r^{2m}\) yields
\begin{equation}\label{eq:Dirac-split}
E_{r,\mathcal D}^0(u)
=
E_{r,\mathcal D}^0(u_\beta)
+
E_{r,\mathcal D}^0(u_\perp).
\end{equation}
This identity requires no bidegree comparison; it uses only the
unitarity of \(U(\varkappa)\) and the \(L^2(S^{\nu-1})\)-orthogonality of
the angular spectral splitting.

\medskip
\emph{ Non-resonant Hardy inequality.}
Expand $u$ in angular eigenmodes: write $u=\sum_\mu \varrho_\mu(r)\,Y_\mu(\varkappa)$, where
$Y_\mu$ runs over an $L^2(S^{\nu-1})$-orthonormal basis of eigensections of
$\mathcal L$ with multiplicity of the eigenvalue $\mu$, and each radial coefficient
$\varrho_\mu\in C^\infty([0,R))$ vanishes for $r$ near $R$ (compact support of $u$) and
is smooth at $r=0$. For such a coefficient \(\varrho=\varrho_\mu\), the cross term in \(\bigl|\varrho'+\tfrac{\mu}{r}\varrho\bigr|^2\) is \(2\frac{\mu}{r}\operatorname{Re}(\varrho'\overline{\varrho})\). Multiplying this cross term by \(r^{2\beta+1}\) and integrating over \([0,R]\), then integrating by parts, gives
$-2\mu\beta\int_0^R r^{2\beta-1}|\varrho|^2\,dr$,
with vanishing boundary terms: the boundary term at \(r=R\) vanishes because \(\varrho\) is compactly supported in \(r<R\), and the boundary term at \(r=0\) vanishes because \(r^{2\beta}|\varrho|^2\to0\) as \(r\to0\) (since \(\varrho\) is bounded near \(0\) and \(\beta>0\)). The same computation with $\mu$ replaced by $\beta$ gives the cross
term of the right-hand side, and subtracting yields
\begin{equation}
\label{eq:sae}
\int_0^R r^{2\beta+1}\Bigl|\varrho'+\tfrac{\mu}{r}\varrho\Bigr|^2dr
=\int_0^R r^{2\beta+1}\Bigl|\varrho'+\tfrac{\beta}{r}\varrho\Bigr|^2dr
+(\mu-\beta)^2\int_0^R r^{2\beta-1}|\varrho|^2dr .
\end{equation}
On $u_\perp$ the spectrum avoids $\beta$, so $|\mu-\beta|\ge1$.  By the
polar formula, the unitarity of \(U(\varkappa)\) and Parseval's identity on
each sphere,
\(\int_{S^{\nu-1}}|\mathcal Du_\perp|^2\,d\varkappa
=\frac14\sum_{\mu\ne\beta}\bigl|\varrho_\mu'+\frac{\mu}{r}\varrho_\mu\bigr|^2\).
Together with
the volume identities $r^{2m}dV=r^{2\beta+1}dr\,d\varkappa$ and
$r^{2m-2}dV=r^{2\beta-1}dr\,d\varkappa$ (recall $2\beta+1=2m+\nu-1$),
\eqref{eq:sae} summed over $\mu\ne\beta$ yields
\begin{equation}
\label{eq:Hnr}
\int_{B_R(z_c)}r^{2m-2}|u_\perp|^2\,dV
\le 4\int_{B_R(z_c)}r^{2m}|\mathcal D u_\perp|^2\,dV .
\end{equation}

\medskip
\emph{Logarithmic Hardy inequality on the resonant channel.}
If \(\mathrm{P}_\beta=0\) (in particular, if
\(\beta\notin\operatorname{spec}(\mathcal L)\)), then \(u_\beta=0\)
and \eqref{eq:Hr} is trivial. We may therefore assume \(u_\beta\ne0\)
in the sequel.

For $u_\beta$, the gap in \eqref{eq:sae} vanishes. For a resonant radial
coefficient $\varrho=\varrho_\beta\in C^\infty([0,R))$, smooth at $r=0$ and vanishing for
$r$ near $R$, set $F(t)=r^\beta \varrho(r)$ with $r=Re^{-t}$, $t\in[0,\infty)$.

At the outer endpoint \(r=R\), i.e.\ \(t=0\), the compact support condition
\(\operatorname{supp}u\Subset B_R(z_c)\) ensures that \(u\) vanishes in a
neighborhood of \(\partial B_R(z_c)\). Since the spectral projection \(\mathrm{P}_\beta\)
acts only on the angular variable and preserves radial support, we have \(\operatorname{supp}u_\beta\Subset B_R(z_c)\), and hence
the radial coefficient \(\lim_{r\to R^-}\varrho(r)=0\), which gives \(F(0)=0\). As \(r\to0\),
i.e.\ \(t\to\infty\), the smoothness of
\(\varrho\) at the origin and \(\beta>0\) imply \(F(t)\to0\) and
\(F'(t)\to0\). Indeed, writing
\[
\varrho(r)=\varrho(0)+r\varrho'(0)+O(r^2),
\]
we obtain
\[
F(t)=r^\beta\varrho(r)
=R^\beta e^{-\beta t}\bigl(\varrho(0)+O(e^{-t})\bigr),
\]
and
\[
F'(t)
=-\beta r^\beta\varrho-r^{\beta+1}\varrho'
=-\beta R^\beta e^{-\beta t}\bigl(\varrho(0)+O(e^{-t})\bigr)
-R^{\beta+1}e^{-(\beta+1)t}\bigl(\varrho'(0)+O(e^{-t})\bigr).
\]
Since \(\beta>0\), both \(F(t)\) and \(F'(t)\) tend to \(0\) as
\(t\to\infty\).
These limits follow from the regularity of $\varrho$ at the singular point, not merely from the compact
support of $u$. Since $\varsigma_R(r)=\varsigma_R(Re^{-t})=1+t$,
\[
\int_0^R\frac{r^{2\beta-1}|\varrho|^2}{\varsigma_R^2}\,dr=\int_0^\infty\frac{|F|^2}{(1+t)^2}\,dt,
\qquad
\int_0^R r^{2\beta+1}\Bigl|\varrho'+\tfrac{\beta}{r}\varrho\Bigr|^2dr=\int_0^\infty|F'|^2dt .
\]

Integration by parts and Cauchy--Schwarz give
\[
\int_0^\infty\frac{|F|^2}{(1+t)^2}\,dt
\le
4\int_0^\infty|F'|^2dt .
\]
Indeed, for \(T>0\),
\(\int_0^T\frac{|F|^2}{(1+t)^2}dt
=-\frac{|F(T)|^2}{1+T}+\int_0^T\frac{2\operatorname{Re}(F'\overline F)}{1+t}dt
\le2\bigl(\int_0^T|F'|^2\bigr)^{1/2}\bigl(\int_0^T\frac{|F|^2}{(1+t)^2}\bigr)^{1/2}\),
using \(F(0)=0\); the last integral is finite because \(F\) is bounded.
Dividing and letting \(T\to\infty\) gives the claim.
Summing over the resonant modes and using
\(u_\beta\in\operatorname{ran}\mathrm{P}_\beta\), so that
\(\mathcal L u_\beta=\beta u_\beta\), we obtain
\[
\int_{B_R(z_c)}\frac{r^{2m-2}}{\varsigma_R^2}|u_\beta|^2\,dV
\le
4\int_{B_R(z_c)}r^{2m}
\Bigl|\Bigl(\partial_r+\frac{\beta}{r}\Bigr)u_\beta\Bigr|^2\,dV .
\]
By the polar decomposition
\(\mathcal D=\frac12U(\varkappa)(\partial_r+\frac1r\mathcal L)\)
and the unitarity of \(U(\varkappa)\),
\[
|\mathcal D u_\beta|^2
=
\frac14
\Bigl|\Bigl(\partial_r+\frac{\beta}{r}\Bigr)u_\beta\Bigr|^2,
\]
hence
\begin{equation}
\label{eq:Hr}
\int_{B_R(z_c)}\frac{r^{2m-2}}{\varsigma_R^2}|u_\beta|^2\,dV
\le 16\int_{B_R(z_c)}r^{2m}|\mathcal D u_\beta|^2\,dV .
\end{equation}

\medskip
\emph{ The weighted energy is controlled by the flat energy.}
Define the two channel energies
\[
\mathcal E_r^\perp(u):=
\int_{B_R(z_c)}r^{2m}|\nabla u_\perp|^2
+\beta^2\!\int_{B_R(z_c)}r^{2m-2}|u_\perp|^2,
\]
\[
\mathcal E_r^\beta(u):=
\int_{B_R(z_c)}\frac{r^{2m}}{\varsigma_R^2}|\nabla u_\beta|^2
+\beta^2\!\int_{B_R(z_c)}\frac{r^{2m-2}}{\varsigma_R^2}|u_\beta|^2 ,
\]
and \(\mathcal E_r(u):=\mathcal E_r^\perp(u)+\mathcal E_r^\beta(u)\).

On the non-resonant channel, apply the multiplier inequality \eqref{GD} with $h=r^m$
($|h'|=m\,r^{m-1}$):
\[
\int_{B_R(z_c)}r^{2m}|\nabla u_\perp|^2
\le
2C_{\mathrm{GD}}^2\int_{B_R(z_c)}r^{2m}|\mathcal D u_\perp|^2
+
2C_{\mathrm{GD}}^2m^2\int_{B_R(z_c)}r^{2m-2}|u_\perp|^2,
\]
the last term controlled by \eqref{eq:Hnr}. On the resonant channel, apply
\eqref{GD} with $h=r^m/\varsigma_R$; since
$|h'|\le(m+1)r^{m-1}/\varsigma_R$,
\[
\int_{B_R(z_c)}\frac{r^{2m}}{\varsigma_R^2}|\nabla u_\beta|^2
\le
2C_{\mathrm{GD}}^2\int_{B_R(z_c)}r^{2m}|\mathcal D u_\beta|^2
+
2C_{\mathrm{GD}}^2(m+1)^2\int_{B_R(z_c)}\frac{r^{2m-2}}{\varsigma_R^2}|u_\beta|^2,
\]
the last term controlled by \eqref{eq:Hr} and the first by the flat energy
since $\varsigma_R\ge1$.

Combining these bounds with \eqref{eq:Hnr}, \eqref{eq:Hr},
\(E_{r,\mathcal D}^0(u_\perp)\le E_{r,\mathcal D}^0(u)\) and
\(E_{r,\mathcal D}^0(u_\beta)\le E_{r,\mathcal D}^0(u)\)
(both consequences of \eqref{eq:Dirac-split}), we obtain
 \[
 \mathcal E_r(u)\le C\,E_{r,\mathcal D}^0(u),
 \qquad C=C(m,n).
 \]
Since \(u\in\Lambda^{p,q}\), the bidegree bridge applied to the original
form \(u\) gives
\[
E_{r,\mathcal D}^0(u)\le E_r^f(u).
\]
Within this proof, this is the only place where the bridge is used.  It is not applied to
\(u_\beta\) or \(u_\perp\), which generally are not of fixed bidegree;
the identity \eqref{eq:Dirac-split}
is a statement about the full Dirac energy and requires no bidegree
comparison.
 Thus
\begin{equation}
\label{eq:Ldec}
 \mathcal E_r(u)\le C\,E_r^f(u).
\end{equation}

\medskip
\emph{Passage to $H^{1/2}$ via Lemma~\ref{lem:SW}.}
Apply Lemma~\ref{lem:SW} (with $\nu=2n\ge2$) componentwise to
$g=r^m u=\theta u$, extended by zero outside $B_R(z_c)$; this is legitimate
because $u\in C_c^\infty(B_R(z_c))$ and $m\ge1$ make $r^m u$ a compactly
supported $H^1$ section. Thus
\begin{equation}
\label{eq:apply-SW}
\|r^m u\|_{H^{1/2}(\R^\nu)}^2
\le C_\nu\Bigl(\|r^m u\|_{L^2}^2+\int_{\R^\nu}|x|\,|\nabla(r^m u)|^2\,dx\Bigr),
\end{equation}
By the product rule $|\nabla(r^m u)|^2\le 2\,r^{2m}|\nabla u|^2+2m^2\,r^{2m-2}|u|^2$,
so with $|x|=r$,
\[
\int_{B_R(z_c)}|x|\,|\nabla(r^m u)|^2
\le 2\int_{B_R(z_c)}r^{2m+1}|\nabla u|^2+2m^2\int_{B_R(z_c)}r^{2m-1}|u|^2 .
\]
Split $u=u_\beta+u_\perp$ and use \eqref{eq:grad-split} with the radial
weights $\kappa=r^{2m+1}$ and $\kappa=r^{2m-1}$, respectively. On the non-resonant
channel, $r\le R$ gives $r^{2m+1}\le R\,r^{2m}$ and
$r^{2m-1}\le R\,r^{2m-2}$, and since $m\le\beta$ (because $n\ge1$) we have
$m^2\le\beta^2$, so
\[
\int_{B_R(z_c)}\!\bigl(r^{2m+1}|\nabla u_\perp|^2+m^2 r^{2m-1}|u_\perp|^2\bigr)
\le R\,\mathcal E_r^\perp(u).
\]
On the resonant channel, write
$r^{2m+1}=(r\varsigma_R^2)\,r^{2m}\varsigma_R^{-2}$ and
$r^{2m-1}=(r\varsigma_R^2)\,r^{2m-2}\varsigma_R^{-2}$. A direct calculation gives
$\sup_{0<r\le R}\,r\varsigma_R(r)^2= 4R/e$ (indeed $r\varsigma_R^2\to0$
as $r\to0$), and again $m^2\le\beta^2$,
\[
\int_{B_R(z_c)}\!\bigl(r^{2m+1}|\nabla u_\beta|^2+m^2 r^{2m-1}|u_\beta|^2\bigr)
\le (4R/e)\,\mathcal E_r^\beta(u).
\]
Adding the two channels,
\begin{equation}
\label{eq:wgrad-bound}
\int_{B_R(z_c)}|x|\,|\nabla(r^m u)|^2
\le (8R/e)\,\mathcal E_r(u).
\end{equation}

Combining  \eqref{eq:Ldec}, \eqref{eq:apply-SW}, \eqref{eq:wgrad-bound}, and
$\|r^m u\|_{L^2}^2=M_r(u)$, and using that the full-space norm
dominates the intrinsic norm on $B_R(z_c)$ up to a constant $C(\nu)$
(Remark~\ref{rem:Heps-norms}),
\[
\|r^m u\|_{H^{1/2}(B_R(z_c))}^2
\le C(\nu)\|r^m u\|_{H^{1/2}(\R^\nu)}^2
\le C\bigl(\mathcal E_r(u)+M_r(u)\bigr)
\le C\bigl(E_r^f(u)+M_r(u)\bigr).
\]
Writing out $E_r^f$ and $M_r$ gives
\[
\|r^m u\|_{H^{1/2}(B_R(z_c))}^2
\le C\Bigl(\|r^m \overline{\partial}_fu\|_{L^2(B_R(z_c))}^2+\|r^m \overline{\partial}_f^{\ast}u\|_{L^2(B_R(z_c))}^2+\|r^m u\|_{L^2(B_R(z_c))}^2\Bigr),
\]
with $C=C(n,m,R)$. This is \eqref{eq:M}.
\end{proof}

We now pass from \(C^\infty_{\Lambda^{p,q}}\) to the \(L^2\)-domain. The first
ingredient is the density of compactly supported smooth forms in the relevant
graph norm.

\begin{lemma}[Smooth compactly supported approximation in the graph norm]
\label{lem:smooth-approx-critical}

 Let \(0\le p\le n\), \(1\le q\le n\), and let
 \[
 u\in L^2_{(p,q)}(B_R(z_c))\cap\dom(\overline{\partial}_f)\cap\dom(\overline{\partial}_f^{\ast}),
 \qquad
 \supp u\Subset B_R(z_c).
 \]
Then there exists a sequence
$u_j\in C_{c,\Lambda^{p,q}}^{\infty}(B_R(z_c))$
such that
\[
u_j\to u\quad\text{in }L^2_{(p,q)}(B_R(z_c)),
\]
and 
$\overline\partial u_j\to\overline\partial u$,
$\overline\partial^{\ast} u_j\to\overline\partial^{\ast} u$
in \(L^2\).
Consequently also
 \[
 \overline{\partial}_fu_j\to \overline{\partial}_fu\quad\text{in }L^2_{(p,q+1)}(B_R(z_c))\oplus L^2_{(p+1,q)}(B_R(z_c)),
 \] \[
 \overline{\partial}_f^{\ast}u_j\to \overline{\partial}_f^{\ast}u\quad\text{in }L^2_{(p,q-1)}(B_R(z_c))\oplus L^2_{(p-1,q)}(B_R(z_c)).
 \]
\end{lemma}

\begin{proof}

Since \(\supp u\Subset B_{R}(z_c)\), the approximation is carried out by
interior mollification, so that no boundary regularity of
\(\partial B_{R}(z_c)\) is involved.  This is the compactly supported
analogue of the graph-norm density results in
\cite[Proposition~2.3]{Straube2010}; no boundary-density issues arise here, and for
convenience, we give the short argument in full below, since we need
the conclusion simultaneously for \(\overline\partial\) and \(\overline\partial^{\ast}\)
(and hence for \(\overline{\partial}_f\) and \(\overline{\partial}_f^{\ast}\)). Fix \(d_c:=\operatorname{dist}(\supp u,\partial
B_{R}(z_c))>0\) and let \(\{\eta_j\}_{j\ge1}\) be a standard mollifier with
\(\supp\eta_j\subset B_{1/j}(0)\), \(\eta_j\ge0\), \(\int\eta_j=1\). For
\(1/j<d_c\) set
\[
u_j:=\eta_j\ast u,
\]
the convolution taken componentwise in a fixed global orthonormal frame of
\(\Lambda^{p,q}\). Then \(u_j\in
C_{c,\Lambda^{p,q}}^{\infty}(B_{R}(z_c))\), the bidegree \((p,q)\) being
preserved because convolution acts on each constant-frame coefficient
separately, and \(u_j\to u\) in \(L^2\).

Recall \(\overline{\partial}_f=\overline\partial+\partial f\wedge\) and
\(\overline{\partial}_f^{\ast}=\overline\partial^{\ast}+\iota_{(\partial f)^\sharp}\). On \(B_{R}(z_c)\), since
\(\supp u\Subset B_{R}(z_c)\), the distributional expressions \(\overline{\partial}_fu\in
L^2_{(p,q+1)}(\Omega)\oplus L^2_{(p+1,q)}(\Omega)\) and \(\overline{\partial}_f^{\ast}u\in L^2_{(p,q-1)}(\Omega)\oplus L^2_{(p-1,q)}(\Omega)\) together
with the boundedness of the zero-order terms give \(\overline\partial u,\overline\partial^{\ast} u\in L^2\)
(subtract the bounded summands \(\partial f\wedge u\in\Lambda^{ p+1,q}\) and
\(\iota_{(\partial f)^\sharp}u\in\Lambda^{p-1,q}\), which lie in bidegrees
different from \(\overline\partial u\in\Lambda^{p,q+1}\) and
\(\overline\partial^{\ast} u\in\Lambda^{p,q-1}\)), and \(\overline\partial^{\ast}u=\vartheta u\) by
Lemma~\ref{lem:domains}(c). The constant-coefficient operators
\(\overline\partial\) and \(\vartheta\) commute with convolution, whence
\[
\overline\partial u_j=\eta_j \ast(\overline\partial u)\to\overline\partial u,
\qquad
\overline\partial^{\ast} u_j=\eta_j \ast(\overline\partial^{\ast} u)\to\overline\partial^{\ast} u
\quad\text{in }L^2 .
\]
For the twisted operators, write
\[
\overline{\partial}_fu_j=\overline\partial u_j+\partial f\wedge u_j,
\qquad
\overline{\partial}_f^{\ast}u_j=\overline\partial^{\ast} u_j+\iota_{(\partial f)^\sharp}u_j .
\]
Since the coefficients of \(\partial f\) are smooth and bounded on
\(\overline{B_{R}(z_c)}\) and \(u_j\to u\) in \(L^2\), the zero-order terms
converge in \(L^2\); combined with the convergences for \(\overline\partial,\overline\partial^{\ast}\) this
gives \(\overline{\partial}_fu_j\to \overline{\partial}_fu\) and \(\overline{\partial}_f^{\ast}u_j\to \overline{\partial}_f^{\ast}u\) in \(L^2\). (When
\( p=0\) the term \(\iota_{(\partial f)^\sharp}u_j\) vanishes and
\(\overline{\partial}_f^{\ast}u_j=\overline\partial^{\ast} u_j\); for \( p\ge1\) it is a genuine bounded zero-order
contribution, handled by the same convergence.)
\end{proof}

\begin{remark}\label{rem:smooth-approx-general}
The proof of Lemma~\ref{lem:smooth-approx-critical} uses only that
\(\supp u\) is a compact subset of the open set on which \(u\) is
considered and that \(\overline\partial\), \(\overline\partial^{\ast}\) (acting on
compactly supported forms as the constant-coefficient formal adjoint)
commute with convolution.  It therefore holds verbatim with
\(B_R(z_c)\) replaced by any open set \(G\subset\Omega\): every
\(u\in L^2_{(p,q)}(\Omega)\cap\dom(\overline{\partial}_f)\cap\dom(\overline{\partial}_f^{\ast})\) with
\(\supp u\Subset G\) is the limit of forms
\(u_j\in C^\infty_{c,\Lambda^{p,q}}(G)\) in the graph norms of
\(\overline{\partial}_f\), \(\overline{\partial}_f^{\ast}\), \(\overline\partial\) and \(\overline\partial^{\ast}\).
\end{remark}

\begin{proposition}[Local \(H^{1/2}\) estimate at an isolated critical point: \(L^2\)-realization case]
\label{prop:critical-L2}
Let \(z_c\in Z\), \(m:=m_c\ge1\), \(0<R\le R_c\), \(0\le p\le n\), \(1\le q\le n\), and let
\[
u\in L^2_{(p,q)}(B_R(z_c))\cap\dom(\overline{\partial}_f)\cap\dom(\overline{\partial}_f^{\ast}),
\qquad
\supp u\Subset B_R(z_c).
\]
Since \(\supp u\Subset B_R(z_c)\), \(u\) may be regarded as a
compactly supported form on \(B_R(z_c)\); hence no boundary term arises from
\(\partial B_R(z_c)\) in the weak formulation.
Then
\[
r^m u\in H^{1/2}(B_R(z_c)),
\]
and there exists a constant \(C=C(n,m,R)\), such that
\begin{equation}
\label{eq:critical-L2}
\|r^m u\|_{H^{1/2}(B_R(z_c))}^2
\le
C\Bigl(
\|r^m \overline{\partial}_fu\|_{L^2(B_R(z_c))}^2
+
\|r^m \overline{\partial}_f^{\ast}u\|_{L^2(B_R(z_c))}^2
+
\|r^m u\|_{L^2(B_R(z_c))}^2
\Bigr).
\end{equation}
\end{proposition}

\begin{proof}
By Lemma~\ref{lem:smooth-approx-critical}, choose
$u_j\in C_{c,\Lambda^{p,q}}^{\infty}(B_R(z_c))$
such that
\[
u_j\to u\quad\text{in }L^2,
\qquad
\overline{\partial}_fu_j\to \overline{\partial}_fu\quad\text{in }L^2,
\qquad
\overline{\partial}_f^{\ast}u_j\to \overline{\partial}_f^{\ast}u\quad\text{in }L^2.
\]
Each \(u_j\) is compactly supported in \(B_R(z_c)\), is smooth, and belongs to
\(\dom(\overline{\partial}_f)\cap\dom(\overline{\partial}_f^{\ast})\). Therefore
Proposition~\ref{prop:M-critical} applies to \(u_j\), and gives
\begin{equation}
\label{eq:critical-uj}
\|r^m u_j\|_{H^{1/2}(B_R(z_c))}^2
\le
C\Bigl(
\|r^m \overline{\partial}_fu_j\|_{L^2(B_R(z_c))}^2
+
\|r^m \overline{\partial}_f^{\ast}u_j\|_{L^2(B_R(z_c))}^2
+
\|r^m u_j\|_{L^2(B_R(z_c))}^2
\Bigr).
\end{equation}

We pass to the limit in \eqref{eq:critical-uj}. Since \(r^m\le R^m\) on
\(B_R(z_c)\), the \(L^2\)-convergences established in
Lemma~\ref{lem:smooth-approx-critical} imply
\[
r^m u_j\to r^m u
\quad\text{in }L^2(B_R(z_c)),
\]
and
\[
r^m \overline{\partial}_fu_j\to r^m \overline{\partial}_fu,
\qquad
r^m \overline{\partial}_f^{\ast}u_j\to r^m \overline{\partial}_f^{\ast}u
\]
in \(L^2(B_R(z_c))\). Hence the right-hand side of \eqref{eq:critical-uj}
converges to
\[
C\Bigl(
\|r^m \overline{\partial}_fu\|_{L^2(B_R(z_c))}^2
+
\|r^m \overline{\partial}_f^{\ast}u\|_{L^2(B_R(z_c))}^2
+
\|r^m u\|_{L^2(B_R(z_c))}^2
\Bigr).
\]

For the left-hand side, \eqref{eq:critical-uj} shows that the sequence
\(\{r^m u_j\}\) is bounded in \(H^{1/2}(B_R(z_c))\). Since \(r^m u_j\to r^m u\)
in \(L^2(B_R(z_c))\), it follows that
\[
r^m u_j\rightharpoonup r^m u
\quad\text{weakly in }H^{1/2}(B_R(z_c))
\]
after passing to a subsequence. The weak limit is \(r^m u\), because the
\(L^2\)-limit is \(r^m u\). By lower semicontinuity of the norm under weak
convergence,
\[
\|r^m u\|_{H^{1/2}(B_R(z_c))}^2
\le
\liminf_{j\to\infty}
\|r^m u_j\|_{H^{1/2}(B_R(z_c))}^2.
\]
Combining this with \eqref{eq:critical-uj} and the convergence of the
right-hand side yields
\[
\|r^m u\|_{H^{1/2}(B_R(z_c))}^2
\le
C\Bigl(
\|r^m \overline{\partial}_fu\|_{L^2(B_R(z_c))}^2
+
\|r^m \overline{\partial}_f^{\ast}u\|_{L^2(B_R(z_c))}^2
+
\|r^m u\|_{L^2(B_R(z_c))}^2
\Bigr).
\]
This is \eqref{eq:critical-L2}. In particular, the finiteness of the
right-hand side implies \(r^m u\in H^{1/2}(B_R(z_c))\). The proof is complete.
\end{proof}

\section{Global \texorpdfstring{\(H^{\epsilon}\)}{H epsilon} Subelliptic Estimate for Fixed-Bidegree Forms}

The following theorems are the global subelliptic statements; their conclusions
are gain of \(\epsilon\) derivatives, uniformly in \( p\).

\begin{theorem}[Global subelliptic estimate]
\label{them:global}
Let \(\Omega\Subset\mathbb C^n\) be a bounded Lipschitz pseudoconvex domain
and let \(f\in\mathcal O(U_\rho)\) and \(Z\) be its isolated interior
critical set.  All assumptions (i)--(iv) hold, and let \(\epsilon\in(0,1/2)\) be the uniform
exponent fixed after \eqref{eq:straube-local-est}. Then there is a constant \(C\), independent of
\(p\), \(q\) and \(u\), such that for every
\(0\le p\le n\), \(1\le q\le n\), and
$u\in L^2_{(p,q)}(\Omega)\cap\dom(\overline{\partial}_f)
  \cap\dom(\overline{\partial}_f^{\ast})$,
we have
\begin{equation}\label{eq:global-unweighted-subelliptic}
\|u\|_{H^{\epsilon}(\Omega)}^2
\le C\Bigl(
  \|\overline{\partial}_f u\|_{L^2(\Omega)}^2
  +\|\overline{\partial}_f^{\ast}u\|_{L^2(\Omega)}^2
  +\|u\|_{L^2(\Omega)}^2
\Bigr).
\end{equation}
The constant depends on the domain and its boundary estimate data, but not
on the bidegree.
\end{theorem}
\begin{proof}
Let \(U_1^{\mathrm{est}},\dots,U_M^{\mathrm{est}}\)
and \(C_1,\dots,C_M\) be the Straube estimate neighborhoods and constants
fixed above (the \(C_\ell\) being the maxima over \(1\le q\le n\)).  The
sets \(\Omega,U_1^{\mathrm{est}},\dots,U_M^{\mathrm{est}}\) form a finite
open cover of the compact set \(\overline\Omega\); choose real functions
\(\psi_0\in C^\infty_c(\Omega)\), \(\psi_\ell\in C^\infty_c(U_\ell^{\mathrm{est}})\),
\(0\le\psi_j\le1\), with \(\psi_0+\sum_{\ell=1}^M\psi_\ell=1\) on
\(\overline\Omega\).  These functions depend only on the geometric data.
Put \(G_\psi:=\max_{0\le j\le M}\|\nabla\psi_j\|_{L^\infty}\).

Let \(u\) be as in the statement.  By Lemma~\ref{lem:domains}(a),(b),
\(u\in\dom(\overline\partial_{(p,q)})\cap\dom(\overline\partial^{\ast}_{(p,q)})\),
and by Lemma~\ref{lem:multiplier-closure-L2} (with \(\sigma=\psi_j|_\Omega\))
together with Lemma~\ref{lem:domains}(a), each \(\psi_ju\) lies in
\(L^2_{(p,q)}(\Omega)\cap\dom(\overline\partial)\cap\dom(\overline\partial^{\ast})\), with
\[
\overline\partial(\psi_ju)=\psi_j\overline\partial u+\overline\partial\psi_j\wedge u,
\qquad
\overline\partial^{\ast}(\psi_ju)=\psi_j\overline\partial^{\ast}u-\iota_{(\overline\partial\psi_j)^\sharp}u
\]
(subtract \(\partial f\wedge\psi_ju\) and \(\iota_{(\partial f)^\sharp}\psi_ju\)
from \eqref{eq:weak-product}).  Since
\(|\overline\partial\psi_j\wedge u|+|\iota_{(\overline\partial\psi_j)^\sharp}u|\le2|\nabla\psi_j||u|\)
pointwise,
\begin{equation}\label{eq:unw-commutator}
\|\overline\partial(\psi_ju)\|+\|\overline\partial^{\ast}(\psi_ju)\|
\le\|\overline\partial u\|+\|\overline\partial^{\ast}u\|+2G_\psi\|u\| .
\end{equation}

\emph{Boundary pieces.}  Fix \(1\le\ell\le M\).  By the consequence part of
Lemma~\ref{lem:componentwise} applied to \eqref{eq:straube-local-est} (with
the uniform exponent \(\epsilon\)), and by \eqref{eq:unw-commutator},
\[
\|\psi_\ell u\|_{H^\epsilon(U_\ell^{\mathrm{est}}\cap\Omega)}
\le\sqrt2\,C_\ell\bigl(\|\overline\partial u\|+\|\overline\partial^{\ast}u\|+2G_\psi\|u\|\bigr),
\]
where the constant \(C_\ell\) already incorporates the quotient-to-intrinsic norm conversion from Remark~\ref{rem:Heps-norms}. In particular, the left-hand side is finite.
Since
\(\supp(\psi_\ell u)\subset\supp\psi_\ell\cap\Omega\) and
\(\operatorname{dist}(\supp\psi_\ell,\\ \mathbb C^n\setminus U_\ell^{\mathrm{est}})>0\),
Lemma~\ref{lem:relative-localization} (with \(U=U_\ell^{\mathrm{est}}\cap\Omega\),
\(g=\psi_\ell u\)) gives
\(\|\psi_\ell u\|_{H^\epsilon(\Omega)}\le C^{\mathrm{loc}}_\ell\|\psi_\ell u\|_{H^\epsilon(U_\ell^{\mathrm{est}}\cap\Omega)}\),
hence
\[
\|\psi_\ell u\|_{H^\epsilon(\Omega)}^2
\le 6\,(C^{\mathrm{loc}}_\ell C_\ell)^2\bigl(1+4G_\psi^2\bigr)
\bigl(\|\overline\partial u\|^2+\|\overline\partial^{\ast}u\|^2+\|u\|^2\bigr).
\]

\emph{Interior piece.}  The form \(v:=\psi_0u\) has compact support in
\(\Omega\).  By Remark~\ref{rem:smooth-approx-general} (with \(G=\Omega\))
there are \(\varphi_j\in C^\infty_{c,\Lambda^{p,q}}(\Omega)\) with
\(\varphi_j\to v\), \(\overline\partial\varphi_j\to\overline\partial v\) and
\(\overline\partial^{\ast}\varphi_j\to\overline\partial^{\ast}v\) in \(L^2\).  By the flat
identity of the inner-product convention, applied to the fixed-bidegree
compactly supported smooth forms \(\varphi_{j}-\varphi_{k}\),
\[
\|\nabla(\varphi_{j}-\varphi_{k})\|^2
=4\|\overline\partial(\varphi_{j}-\varphi_{k})\|^2+4\|\overline\partial^{\ast}(\varphi_{j}-\varphi_{k})\|^2
\longrightarrow0 ,\quad \text{for all}\,\, j,k\rightarrow \infty,
\]
so \((\varphi_j)\) is Cauchy in \(H^1(\mathbb C^n)\) (zero extensions); its
limit is \(v\), whence \(v\in H^1(\mathbb C^n)\) and
\(\|\nabla v\|^2=4\|\overline\partial v\|^2+4\|\overline\partial^{\ast}v\|^2\).
By \eqref{eq:H1-Heps} and \eqref{eq:unw-commutator},
\[
\|v\|_{H^\epsilon(\Omega)}^2
\le(1+c_{2n,\epsilon})\bigl(\|v\|^2+4\|\overline\partial v\|^2+4\|\overline\partial^{\ast}v\|^2\bigr)
\le 13(1+c_{2n,\epsilon})(1+4G_\psi^2)
\bigl(\|\overline\partial u\|^2+\|\overline\partial^{\ast}u\|^2+\|u\|^2\bigr).
\]

Since \(u=\psi_0u+\sum_{\ell=1}^M\psi_\ell u\)
on \(\Omega\) and the Gagliardo norm satisfies the triangle inequality,
\(\|u\|_{H^\epsilon(\Omega)}^2\le(M+1)\sum_{j=0}^M\|\psi_ju\|_{H^\epsilon(\Omega)}^2\).
Inserting the two bounds above and then the bridge \eqref{eq:bridge},
\(\|\overline\partial u\|^2+\|\overline\partial^{\ast}u\|^2\le\|\overline{\partial}_fu\|^2+\|\overline{\partial}_f^{\ast}u\|^2\)
(valid because \(u\) has fixed bidegree \((p,q)\)), we obtain
\eqref{eq:global-unweighted-subelliptic} with a constant depending only
on \(n\), \(\epsilon\), \(M\), the \(C_\ell\), the \(C^{\mathrm{loc}}_\ell\) and
\(G_\psi\), hence not on \(p\), \(q\) or \(u\).  The critical points lie in
the interior piece and cause no loss, since the twist is of order zero and
is discarded through the bridge.
\end{proof}

\begin{theorem}[Weighted global subelliptic estimate]
\label{thm:global-Hs}
Let \(\Omega\Subset\mathbb C^n\) be a bounded Lipschitz pseudoconvex domain
and let \(f\in\mathcal O(U_\rho)\) and \(Z\) be its isolated interior
critical set.  All assumptions (i)--(iv) hold, and let \(\epsilon\in(0,1/2)\) be the uniform
exponent fixed after \eqref{eq:straube-local-est}.
Let \(\theta\in W^{1,\infty}(\Omega)\) be the order-adapted weight.  Set
\begin{equation}\label{eq:theta-properties}
d_1:=\min_{\overline\Omega\setminus\bigcup_a B_{R_a/2}(z_a)}\theta>0,
\qquad
D_{\theta}:=\|\theta\|_{L^\infty(\Omega)}<\infty .
\end{equation}
Then there exists a constant \(C\), depending only on \(\Omega\), the
geometric data, and \(\theta\), but not on \(p\), \(q\) or \(u\), such that
for all \(0\le p\le n\), \(1\le q\le n\) and all
$
u\in L^2_{(p,q)}(\Omega)\cap \dom(\overline{\partial}_f)\cap
\dom(\overline{\partial}_f^{\ast})$
(that is, all \((p,q)\)-forms in \(\dom(\overline{\partial}_f)\) satisfying the
twisted \(\overline{\partial}_f\)-Neumann boundary condition in the weak sense),
\begin{equation}
\label{eq:global-subelliptic}
\|\theta u\|_{H^{\epsilon}(\Omega)}^2
\le C\Bigl(
\|\theta \overline{\partial}_f u\|_{L^2(\Omega)}^2
+\|\theta \overline{\partial}_f^{\ast} u\|_{L^2(\Omega)}^2
+\|\theta u\|_{L^2(\Omega)}^2
\Bigr).
\end{equation}
\end{theorem}

\begin{proof}
\emph{ The boundary collar.}
By the construction of the Straube estimate neighborhoods, there exist finitely many
boundary points \(P_1,\dots,P_M\in\partial\Omega\), estimate
neighborhoods \(U_\ell^{\mathrm{est}}\Subset U_{P_\ell}\) of
\(P_\ell\), and constants \(C_\ell>0\) such that
\[
\partial\Omega\subset\bigcup_{\ell=1}^M U_\ell^{\mathrm{est}},
\]
and the local \((0,q)\) estimate with uniform exponent \(\epsilon\)
holds on each \(U_\ell^{\mathrm{est}}\cap\Omega\) for the class in
\eqref{eq:local-estimate-Uj-prop}.
Since
\(\bigcup_a\overline{B_{R_a}(z_a)}\Subset\Omega\), there is
\(\tau_\partial>0\) such that
\[
\bigl\{z:\operatorname{dist}(z,\partial\Omega)<\tau_\partial\bigr\}
\cap
\bigcup_{a=1}^N\overline{B_{R_a}(z_a)}=\varnothing .
\]
By the shrinking lemma applied to the finite open cover
\(\{U_\ell^{\mathrm{est}}\}\) of \(\partial\Omega\), and then
intersecting with the thin collar
\(\{\operatorname{dist}(\cdot,\partial\Omega)<\tau_\partial/2\}\), we
obtain open sets \(V_\ell\) and, after relabeling, estimate
neighborhoods \(U_\ell^{\mathrm{est}}\) such that
\[
V_\ell\cap\partial\Omega\ne\varnothing,\qquad
\overline{V_\ell}\subset U_\ell^{\mathrm{est}},\qquad
\partial\Omega\subset\bigcup_{\ell=1}^M V_\ell,
\]
and
\[
\overline{V_\ell}\cap
\bigcup_{a=1}^N\overline{B_{R_a}(z_a)}=\varnothing .
\]
Put
\[
V_\partial:=\bigcup_{\ell=1}^M V_\ell .
\]
Then
\[
\overline{V_\partial}\cap
\bigcup_{a=1}^N\overline{B_{R_a}(z_a)}=\varnothing,
\]
and therefore, by \eqref{eq:theta-properties},
$0<d_1\le\theta\le D_{\theta}$
on $V_\ell\cap\overline\Omega$
for every \(\ell\).  Since \(\overline{V_\ell}\cap Z=\varnothing\) and
\(\overline{V_\ell}\subset U_\ell^{\mathrm{est}}\),
Proposition~\ref{prop:subelliptic-collar-piece-L2} applies with
\(V=V_\ell\), \(K=1\) and \(U_1=U_\ell^{\mathrm{est}}\).  It gives an open set
\(W_\ell\subset V_\ell\) with
\[
V_\ell\cap\partial\Omega\subset W_\ell,
\]
independent of the form, and a constant \(\widetilde C_\ell\) (not to be
confused with the Straube constant \(C_\ell\)), independent of the
form, such that for every \((p,q)\)-form
\(v\in L^2_{(p,q)}(\Omega)\cap\dom(\overline{\partial}_f)
\cap\dom(\overline{\partial}_f^{\ast})\) satisfying the weak
\(\overline{\partial}_f\)-Neumann boundary condition,
\begin{equation}\label{eq:collar-general}
\|\theta v\|_{H^\epsilon(W_\ell\cap\Omega)}^2
\le \widetilde C_\ell\Bigl(
\|\overline{\partial}_f(\theta v)\|_{L^2(\Omega)}^2
+\|\overline{\partial}_f^{\ast}(\theta v)\|_{L^2(\Omega)}^2
+\|\theta v\|_{L^2(\Omega)}^2
\Bigr).
\end{equation}
Since \(V_\ell\cap\partial\Omega\subset W_\ell\), the family
\(\{W_\ell\}\) still covers \(\partial\Omega\).  Put
$W_\partial:=\bigcup_{\ell=1}^M W_\ell$.
Because \(W_\ell\subset V_\ell\) and
\(\overline{V_\partial}\) is disjoint from
\(\bigcup_a\overline{B_{R_a}(z_a)}\), we have
\[
\overline{W_\partial}\cap
\bigcup_{a=1}^N\overline{B_{R_a}(z_a)}=\varnothing.
\]
Therefore, there exists a compact neighborhood \(K_\partial\) of
\(\partial\Omega\) with
$\partial\Omega\subset K_\partial\subset W_\partial$.

\medskip
\emph{ The partition of unity.}
Set
\[
U_0:=\Omega\setminus\Bigl(
\bigcup_{a=1}^N\overline{B_{R_a/2}(z_a)}\cup K_\partial\Bigr).
\]
Then
$\{U_0,B_{R_1}(z_1),\dots,B_{R_N}(z_N),W_\partial\}$
is an open cover of \(\overline\Omega\) in the relative topology.
Choose a smooth partition of unity
\[
\chi_0+\sum_{a=1}^N\chi_a+\chi_\partial\equiv1
\quad\text{on }\overline\Omega
\]
subordinate to this cover, with
\[
\chi_a\in C_c^\infty(B_{R_a}(z_a)),\qquad
\chi_a\equiv1\quad\text{on }B_{R_a/2}(z_a),
\]
\[
\chi_\partial\in C_c^\infty(W_\partial),\qquad
\chi_0\in C_c^\infty(U_0).
\]
Such a partition exists: first choose nonnegative bumps \(\tilde\chi_j\)
with the required supports, whose positivity sets
\(\{\tilde\chi_j>0\}\) cover the compact set \(\overline\Omega\)
(shrinking lemma for the finite open cover above), and with \(\tilde\chi_a\equiv1\) on
\(B_{R_a/2}(z_a)\), then normalize by
\(S=\tilde\chi_0+\sum_a\tilde\chi_a+\tilde\chi_\partial>0\).  Since
\(U_0\) and \(W_\partial\) avoid \(B_{R_a/2}(z_a)\), and
\(\operatorname{supp}\tilde\chi_b\subset B_{R_b}(z_b)\) is disjoint from
\(B_{R_a}(z_a)\) for \(b\ne a\) (the balls \(B_{2R_a}(z_a)\) being pairwise
disjoint), we have \(S=\tilde\chi_a=1\) on
\(B_{R_a/2}(z_a)\), hence \(\chi_a=1\) there.

Observe that
\(\operatorname{supp}\nabla\chi_a\subset
B_{R_a}(z_a)\setminus B_{R_a/2}(z_a)\), and on this annulus
\[
\theta=r_a^{m_a}\ge (R_a/2)^{m_a}\ge d_1.
\]
Moreover
\(\operatorname{supp}\nabla\chi_0\subset U_0
\subset\Omega\setminus\bigcup_aB_{R_a/2}(z_a)\),
so \(\theta\ge d_1\) there, and
\(\operatorname{supp}\nabla\chi_\partial\subset W_\partial\),
where \(\theta\ge d_1\) as well.  Hence, for every index
\(j\in\{0,\partial,1,\dots,N\}\),
$\theta\ge d_1>0$
on $\operatorname{supp}\nabla\chi_j$.
Consequently there is a constant \(C_\chi\), depending on \(\theta\) and
the cutoffs, such that
\begin{equation}\label{eq:cutoff-commutator-bound}
\|\theta(\overline\partial\chi_j)\wedge u\|_{L^2}
+\|\theta\iota_{(\overline{\partial}\chi_j)^\sharp}u\|_{L^2}
\le C_\chi\|\theta u\|_{L^2}
\end{equation}
for \(j\in\{0,\partial,1,\dots,N\}\).  Indeed,
\[
\|\theta(\overline\partial\chi_j)\wedge u\|_{L^2}
\le\|\theta\|_{L^\infty}\|\nabla\chi_j\|_{L^\infty}
\|u\|_{L^2(\operatorname{supp}\nabla\chi_j)}
\le D_{\theta}\|\nabla\chi_j\|_{L^\infty}d_1^{-1}\|\theta u\|_{L^2},
\]
and similarly for the contraction term; hence one may take
$C_\chi:=2D_{\theta}\Bigl(\max_j\|\nabla\chi_j\|_{L^\infty}\Bigr)d_1^{-1}$.

\medskip
\emph{ The boundary collar estimate.}
Choose a smooth partition of unity
\[
\psi_\ell\in C_c^\infty(W_\ell),\quad 0\le\psi_\ell\le1,
\qquad
\sum_{\ell=1}^M\psi_\ell=1
\quad\text{on }\operatorname{supp}\chi_\partial.
\]
Set
$w_\ell:=\psi_\ell\chi_\partial u$.
Then
$\operatorname{supp}w_\ell\subset W_\ell\cap\Omega\subset V_\ell\cap\Omega$.
Since \(u\) is a \((p,q)\)-form, so is \(w_\ell\), and
$\theta\chi_\partial u=\sum_{\ell=1}^M\theta w_\ell$.
By Lemma~\ref{lem:multiplier-closure-L2},
\(w_\ell\in\dom(\overline{\partial}_f)\cap\dom(\overline{\partial}_f^{\ast})\), and \(w_\ell\) satisfies the
twisted \(\overline{\partial}_f\)-Neumann boundary condition.  Hence
\(\theta w_\ell\in\dom(\overline{\partial}_f)\cap\dom(\overline{\partial}_f^{\ast})\).  Since
\(0<d_1\le\theta\le D_{\theta}\) on \(V_\ell\cap\overline\Omega\), and
\(\operatorname{supp}w_\ell\subset V_\ell\cap\Omega\), we may apply the
general collar estimate \eqref{eq:collar-general} with \(v=w_\ell\).
This gives
\begin{equation}\label{eq:boundary-collar}
\|\theta w_\ell\|_{H^\epsilon(W_\ell\cap\Omega)}^2
\le \widetilde C_\ell\Bigl(
\|\overline{\partial}_f(\theta w_\ell)\|_{L^2(\Omega)}^2
+\|\overline{\partial}_f^{\ast}(\theta w_\ell)\|_{L^2(\Omega)}^2
+\|\theta w_\ell\|_{L^2(\Omega)}^2
\Bigr).
\end{equation}
By the weak product rule \eqref{eq:weak-product},
\[
\overline{\partial}_f(\theta w_\ell)
=\theta \overline{\partial}_fw_\ell+(\overline\partial\theta)\wedge w_\ell,
\qquad
\overline{\partial}_f^{\ast}(\theta w_\ell)
=\theta \overline{\partial}_f^{\ast}w_\ell-\iota_{(\overline{\partial}\theta)^\sharp}w_\ell.
\]
Since \(w_\ell\) is supported in \(W_\ell\), on which
\(\theta\ge d_1>0\), we have
\[
\|(\overline\partial\theta)\wedge w_\ell\|_{L^2}
\le\|\nabla\theta\|_{L^\infty}\|w_\ell\|_{L^2}
\le d_1^{-1}\|\nabla\theta\|_{L^\infty}\|\theta w_\ell\|_{L^2},
\]
and similarly for the contraction.  Hence
\[
\|\overline{\partial}_f(\theta w_\ell)\|_{L^2}
\le\|\theta \overline{\partial}_fw_\ell\|_{L^2}+C\|\theta w_\ell\|_{L^2},
\]
and
\[
\|\overline{\partial}_f^{\ast}(\theta w_\ell)\|_{L^2}
\le\|\theta \overline{\partial}_f^{\ast}w_\ell\|_{L^2}+C\|\theta w_\ell\|_{L^2}.
\]
Now
\[
\overline{\partial}_fw_\ell=\psi_\ell\chi_\partial \overline{\partial}_fu+
(\overline\partial(\psi_\ell\chi_\partial))\wedge u,
\]
and similarly for \(\overline{\partial}_f^{\ast}\).  The commutator terms are supported in
\(\operatorname{supp}\psi_\ell\subset W_\ell\), where \(\theta\ge d_1>0\);
hence, by the argument proving \eqref{eq:cutoff-commutator-bound} with
\(\chi_j\) replaced by \(\psi_\ell\chi_\partial\),
\[
\|\theta(\overline\partial(\psi_\ell\chi_\partial))\wedge u\|_{L^2}
+\|\theta\iota_{(\overline{\partial}(\psi_\ell\chi_\partial))^\sharp}u\|_{L^2}
\le D_{\theta}\,\|\nabla(\psi_\ell\chi_\partial)\|_{L^\infty}\,d_1^{-1}\|\theta u\|_{L^2}.
\]
Together with \(0\le\psi_\ell\chi_\partial\le1\), this gives
\[
\|\theta \overline{\partial}_fw_\ell\|_{L^2}
\le\|\theta \overline{\partial}_fu\|_{L^2}+C\|\theta u\|_{L^2},
\]
and
\[
\|\theta \overline{\partial}_f^{\ast}w_\ell\|_{L^2}
\le\|\theta \overline{\partial}_f^{\ast}u\|_{L^2}+C\|\theta u\|_{L^2}.
\]
Substituting these bounds into \eqref{eq:boundary-collar} yields, for
each \(\ell\),
\begin{equation}\label{eq:boundary-local-final}
\|\theta w_\ell\|_{H^\epsilon(W_\ell\cap\Omega)}^2
\le \widetilde C_\ell'\Bigl(
\|\theta \overline{\partial}_fu\|_{L^2(\Omega)}^2
+\|\theta \overline{\partial}_f^{\ast}u\|_{L^2(\Omega)}^2
+\|\theta u\|_{L^2(\Omega)}^2
\Bigr).
\end{equation}
Since \(\psi_\ell\in C_c^\infty(W_\ell)\), we have
$d_\ell:=\operatorname{dist}\bigl(\operatorname{supp}\psi_\ell,
\partial W_\ell\bigr)>0$.
Thus
$\operatorname{supp}(\theta w_\ell)
\subset\operatorname{supp}\psi_\ell\cap\Omega$,
and consequently
$\operatorname{dist}\bigl(\operatorname{supp}(\theta w_\ell),
\Omega\setminus W_\ell\bigr)\ge d_\ell>0$.

Applying the relative localization estimate of
Lemma~\ref{lem:relative-localization} with
\(U=W_\ell\cap\Omega\) and \(g=\theta w_\ell\), we obtain
\begin{equation}\label{eq:boundary-zero-extension}
\|\theta w_\ell\|_{H^\epsilon(\Omega)}^2
\le C_{\mathrm{loc},\ell}
\|\theta w_\ell\|_{H^\epsilon(W_\ell\cap\Omega)}^2.
\end{equation}
Combining \eqref{eq:boundary-local-final} and
\eqref{eq:boundary-zero-extension}, and using
\[
\theta\chi_\partial u=\sum_{\ell=1}^M\theta w_\ell,
\]
together with the triangle inequality for the \(H^\epsilon\)-norm, we
obtain
\begin{equation}\label{eq:boundary-final}
\|\theta\chi_\partial u\|_{H^\epsilon(\Omega)}^2
\le M\sum_{\ell=1}^M\|\theta w_\ell\|_{H^\epsilon(\Omega)}^2
\le C\Bigl(
\|\theta \overline{\partial}_fu\|_{L^2(\Omega)}^2
+\|\theta \overline{\partial}_f^{\ast}u\|_{L^2(\Omega)}^2
+\|\theta u\|_{L^2(\Omega)}^2
\Bigr).
\end{equation}

\medskip
\emph{ Interior isolated zeros.}
Fix \(z_a\in Z\), so \(m_a\ge1\). 
Since \(\chi_a\in C_c^\infty(B_{R_a}(z_a))\subset W^{1,\infty}\)
and \(\supp\chi_a\Subset B_{R_a}(z_a)\Subset\Omega\),
Lemma~\ref{lem:multiplier-closure-L2} gives
\(v_a=\chi_a u\in L^2_{(p,q)}\cap\dom(\overline{\partial}_f)\cap\dom(\overline{\partial}_f^{\ast})\)
with \(\supp v_a\Subset B_{R_a}(z_a)\); by Lemma~\ref{lem:domains}(c),
\(v_a\) satisfies the hypotheses of Proposition~\ref{prop:critical-L2}
on \(B_{R_a}(z_a)\), where \(\theta=r_a^{m_a}\).
By
Proposition~\ref{prop:critical-L2}, applied with \(m=m_a\) to the fixed
\((p,q)\)-form \(v_a\),
\begin{equation}\label{eq:critical-va}
\|\theta v_a\|_{H^{1/2}(B_{R_a})}^2
\le C_a\Bigl(
\|\theta \overline{\partial}_fv_a\|_{L^2}^2
+\|\theta \overline{\partial}_f^{\ast}v_a\|_{L^2}^2
+\|\theta v_a\|_{L^2}^2
\Bigr).
\end{equation}
We estimate the three terms on the right-hand side.  Since
\[
\overline{\partial}_fv_a=\chi_a\overline{\partial}_fu+(\overline\partial\chi_a)\wedge u,
\qquad
\overline{\partial}_f^{\ast}v_a=\chi_a\overline{\partial}_f^{\ast}u-\iota_{(\overline{\partial}\chi_a)^\sharp}u,
\]
we have
\[
\|\theta \overline{\partial}_fv_a\|_{L^2}
\le\|\theta \overline{\partial}_fu\|_{L^2}
+\|\theta(\overline\partial\chi_a)\wedge u\|_{L^2},
\]
and
\[
\|\theta \overline{\partial}_f^{\ast}v_a\|_{L^2}
\le\|\theta \overline{\partial}_f^{\ast}u\|_{L^2}
+\|\theta\iota_{(\overline{\partial}\chi_a)^\sharp}u\|_{L^2}.
\]
The commutator terms are controlled by
\eqref{eq:cutoff-commutator-bound}:
\[
\|\theta(\overline\partial\chi_a)\wedge u\|_{L^2}
+\|\theta\iota_{(\overline{\partial}\chi_a)^\sharp}u\|_{L^2}
\le C_\chi\|\theta u\|_{L^2}.
\]
Moreover
$\|\theta v_a\|_{L^2}\le\|\theta u\|_{L^2}$.
Substituting these bounds into \eqref{eq:critical-va} yields
\begin{equation}\label{eq:critical-H12}
\|\theta\chi_au\|_{H^{1/2}(B_{R_a})}^2
\le C_a\Bigl(
\|\theta \overline{\partial}_fu\|_{L^2(\Omega)}^2
+\|\theta \overline{\partial}_f^{\ast}u\|_{L^2(\Omega)}^2
+\|\theta u\|_{L^2(\Omega)}^2
\Bigr).
\end{equation}
Since \(\operatorname{supp}\chi_a\Subset B_{R_a}\) and
\(\epsilon<1/2\), the inequality \eqref{eq:Heps-monotone} with
\(\epsilon^{\prime}=1/2\), together with the
relative-localization estimate of Lemma~\ref{lem:relative-localization}
(applied with \(U=B_{R_a}\) and \(g=\theta\chi_a u\)) gives
\[
\|\theta\chi_au\|_{H^\epsilon(\Omega)}
\le C_a\|\theta\chi_au\|_{H^\epsilon(B_{R_a})}
\le C_a'\|\theta\chi_au\|_{H^{1/2}(B_{R_a})}.
\]
Combining this with \eqref{eq:critical-H12}, we obtain
\begin{equation}\label{eq:critical-final}
\|\theta\chi_au\|_{H^\epsilon(\Omega)}^2
\le C\Bigl(
\|\theta \overline{\partial}_fu\|_{L^2(\Omega)}^2
+\|\theta \overline{\partial}_f^{\ast}u\|_{L^2(\Omega)}^2
+\|\theta u\|_{L^2(\Omega)}^2
\Bigr).
\end{equation}
\medskip
\emph{ The interior region away from \(Z\).}
On \(\operatorname{supp}\chi_0\) the weight \(\theta\) satisfies
$d_1\le\theta\le D_{\theta}$.
Since \(\chi_0\in C_c^\infty(U_0)\subset W^{1,\infty}\),
Lemma~\ref{lem:multiplier-closure-L2} gives
\(\chi_0u\in L^2_{(p,q)}\cap\dom(\overline{\partial}_f)\cap\dom(\overline{\partial}_f^{\ast})\), and the domain
equality \eqref{eq:domain-equality}
gives \(\chi_0u\in\dom(\overline\partial)\cap\dom(\overline\partial^{\ast})\).
Moreover, \(\chi_0u\) is compactly supported in \(U_0\subset\Omega\).

Exactly as for the interior piece \(\psi_0u\) in the proof of
Theorem~\ref{them:global} (Remark~\ref{rem:smooth-approx-general} with
\(G=U_0\), and the flat identity of the inner-product convention applied
to the Cauchy sequence of approximants), the zero extension of \(\chi_0u\)
lies in \(H^1(\mathbb C^n)\) and
\[
\|\chi_0u\|_{H^1(\Omega)}^2
=
\|\chi_0u\|_{L^2(\Omega)}^2
+4\|\overline\partial(\chi_0u)\|_{L^2(\Omega)}^2
+4\|\overline\partial^{\ast}(\chi_0u)\|_{L^2(\Omega)}^2.
\]
Since \(\chi_0u\) is a fixed \((p,q)\)-form, the bidegree bridge gives
\[
\|\overline\partial(\chi_0u)\|_{L^2}^2+\|\overline\partial^{\ast}(\chi_0u)\|_{L^2}^2
\le
\|\overline{\partial}_f(\chi_0u)\|_{L^2}^2+\|\overline{\partial}_f^{\ast}(\chi_0u)\|_{L^2}^2,
\]
and hence
\[
\|\chi_0u\|_{H^1(\Omega)}^2
\le C_0\Bigl(
\|\overline{\partial}_f(\chi_0u)\|_{L^2(\Omega)}^2
+\|\overline{\partial}_f^{\ast}(\chi_0u)\|_{L^2(\Omega)}^2
+\|\chi_0u\|_{L^2(\Omega)}^2
\Bigr).
\]
Multiplying by the weight \(\theta\), we use the product rule
$\nabla(\theta\chi_0u)=\theta\nabla(\chi_0u)+\chi_0u\nabla\theta$,
so that
\[
\|\theta\chi_0u\|_{H^1(\Omega)}^2
\le 2\|\theta\|_{L^\infty}^2\|\chi_0u\|_{H^1(\Omega)}^2
+2\|\nabla\theta\|_{L^\infty}^2\|\chi_0u\|_{L^2(\Omega)}^2.
\]
Combining this with \(d_1\le\theta\le D_{\theta}\) on
\(\operatorname{supp}\chi_0\), we obtain
\[
\|\theta\chi_0u\|_{H^1(\Omega)}^2
\le C(\|\theta\|_{W^{1,\infty}},d_1,D_{\theta})\Bigl(
\|\theta \overline{\partial}_f(\chi_0u)\|_{L^2(\Omega)}^2
+\|\theta \overline{\partial}_f^{\ast}(\chi_0u)\|_{L^2(\Omega)}^2
+\|\theta\chi_0u\|_{L^2(\Omega)}^2
\Bigr).
\]
Using
\[
\overline{\partial}_f(\chi_0u)=\chi_0\overline{\partial}_fu+(\overline\partial\chi_0)\wedge u,
\qquad
\overline{\partial}_f^{\ast}(\chi_0u)=\chi_0\overline{\partial}_f^{\ast}u-\iota_{(\overline{\partial}\chi_0)^\sharp}u,
\]
and \eqref{eq:cutoff-commutator-bound}, together with the quantitative
embedding \eqref{eq:H1-Heps} of Remark~\ref{rem:Heps-norms} applied to the
zero extension of the compactly supported form \(\theta\chi_0u\in H^1\)
(which gives
\(\|\theta\chi_0u\|_{H^\epsilon(\Omega)}^2\le(1+c_{2n,\epsilon})\|\theta\chi_0u\|_{H^1(\Omega)}^2\)),
we obtain
\begin{equation}\label{eq:interior-final}
\|\theta\chi_0u\|_{H^\epsilon(\Omega)}^2
\le C\Bigl(
\|\theta \overline{\partial}_fu\|_{L^2(\Omega)}^2
+\|\theta \overline{\partial}_f^{\ast}u\|_{L^2(\Omega)}^2
+\|\theta u\|_{L^2(\Omega)}^2
\Bigr).
\end{equation}

\medskip
\emph{ Conclusion.}
Since
\[
\theta u=\theta\chi_0u+\sum_{a=1}^N\theta\chi_au
+\theta\chi_\partial u,
\]
the triangle inequality gives
\(\|\theta u\|_{H^\epsilon(\Omega)}^2\le(N+2)\bigl(\|\theta\chi_0u\|_{H^\epsilon(\Omega)}^2
+\sum_a\|\theta\chi_au\|_{H^\epsilon(\Omega)}^2+\|\theta\chi_\partial u\|_{H^\epsilon(\Omega)}^2\bigr)\).
Here \(\|\cdot\|_{H^\epsilon(\Omega)}\) is the intrinsic Gagliardo norm of
Remark~\ref{rem:Heps-norms}, which satisfies the triangle inequality
because the Gagliardo seminorm is an \(L^2(\Omega\times\Omega)\)-norm of
\((F(x)-F(y))|x-y|^{-n-\epsilon}\).  The constants from each piece
depend only on the geometry and \(\theta\), which are fixed.  Summing
over the finitely many pieces in the partition of unity gives the global
constant \(C\).

Inserting \eqref{eq:boundary-final}, \eqref{eq:critical-final} and
\eqref{eq:interior-final} yields
\[
\|\theta u\|_{H^\epsilon(\Omega)}^2
\le C\Bigl(
\|\theta \overline{\partial}_fu\|_{L^2(\Omega)}^2
+\|\theta \overline{\partial}_f^{\ast}u\|_{L^2(\Omega)}^2
+\|\theta u\|_{L^2(\Omega)}^2
\Bigr).
\]
This is \eqref{eq:global-subelliptic}.  The proof is complete.
\end{proof}

\noindent\textbf{Acknowledgments}

\vspace{0.5em}
\noindent 
This work was partially supported by NSFC 11501432. The author is very grateful to Professor Huijun Fan for many discussions and for providing the reference \cite{WeFan}. We thank Professor Emil J. Straube for helpful discussions \cite{Straube2010}. AI tools were used to proofread the final version of the paper, leading to several minor corrections.

\end{document}